\documentclass[11pt,reqno]{article}
\usepackage{amssymb, amsmath, amsthm, amsfonts, amscd, epsfig}
\usepackage[english]{babel}
\usepackage{graphicx, epsfig}
\usepackage{color}

\newtheorem{theorem}{Theorem}[section]
\newtheorem{cor}[theorem]{Corollary}
\newtheorem{lemma}[theorem]{Lemma}
\newtheorem{prop}[theorem]{Proposition}

\newcommand{\ds}{\displaystyle}
\newcommand{\pf}{\noindent {\sl Proof}. \ }
\newcommand{\p}{\partial}

\newcommand{\eqnref}[1]{(\ref {#1})}

\newcommand{\Rbb}{\mathbb{R}}

\newcommand{\Ecal}{\mathcal{E}}

\newcommand{\Ical}{\mathcal{I}}
\newcommand{\Jcal}{\mathcal{J}}
\newcommand{\Kcal}{\mathcal{K}}

\newcommand{\Mcal}{\mathcal{M}}

\newcommand{\Qcal}{\mathcal{Q}}

\def\Be{{\bf e}}

\def\Bg{{\bf g}}
\def\Bh{{\bf h}}

\def\Bp{{\bf p}}

\def\Bu{{\bf u}}
\def\Bv{{\bf v}}

\def\Bx{{\bf x}}

\def\BA{{\bf A}}

\def\BC{{\bf C}}

\def\BU{{\bf U}}

\newcommand{\Ga}{\alpha}
\newcommand{\Gb}{\beta}
\newcommand{\Gd}{\delta}

\newcommand{\Gg}{\gamma}

\newcommand{\Gm}{\mu}

\newcommand{\Gt}{\theta}

\newcommand{\Gs}{\sigma}

\newcommand{\Gy}{\psi}
\newcommand{\Gz}{\zeta}
\newcommand{\GD}{\Delta}

\newcommand{\GG}{\Gamma}

\newcommand{\GO}{\Omega}

\newcommand{\BGG}{{\bf \GG}}

\newcommand{\Bpsi}{\mbox{\boldmath $\Gy$}}

\newcommand{\beq}{\begin{equation}}
\newcommand{\eeq}{\end{equation}}

\newcommand{\BUst}{\BU_{\rm ex}}
\newcommand{\BUsh}{\BU_{\rm sh}}
\newcommand{\Bust}{\Bu_{\rm ex}}
\newcommand{\Bush}{\Bu_{\rm sh}}
\newcommand{\pst}{p_{\rm ex}}
\newcommand{\psh}{p_{\rm sh}}

\numberwithin{equation}{section}
\numberwithin{figure}{section}

\begin{document}

\title{Quantitative estimates for stress concentration of the Stokes flow between adjacent circular cylinders\thanks{This work is supported by NRF grants No. 2017R1A4A1014735, 2019R1A2B5B01069967 and 2020R1C1C1A01010882.}}

\author{Habib Ammari\thanks{Department of Mathematics, ETH Z\"urich, R\"amistrasse 101, CH-8092 Z\"urich, Switzerland (habib.ammari@math.ethz.ch).} \and Hyeonbae Kang\thanks{Department of Mathematics and Institute of Mathematics, Inha University, Incheon 22212, S. Korea (hbkang@inha.ac.kr, dokim@inha.ac.kr). } \and Do Wan Kim\footnotemark[3] \and Sanghyeon Yu\thanks{Department of Mathematics, Korea University, Seoul 02841, S. Korea (sanghyeon\_yu@korea.ac.kr).}}

\date{}
\maketitle

\begin{abstract}
When two inclusions with high contrast material properties are located close to each other in a homogeneous medium, stress may become arbitrarily large in the narrow region between them. In this paper, we investigate such stress concentration in the two-dimensional Stokes flow when inclusions are the two-dimensional cross sections of circular cylinders of the same radii and the background velocity field is linear. We construct two vector-valued functions which completely capture the singular behavior of the stress and derive an asymptotic representation formula for the stress in terms of these functions as the distance between the two cylinders tends to zero. We then show, using the representation formula, that the stress always blows up by proving that either the pressure or the shear stress component of the stress tensor blows up. The blow-up rate is shown to be $\Gd^{-1/2}$, where $\Gd$ is the distance between the two cylinders. To our best knowledge, this work is the first to rigorously derive the asymptotic solution in the narrow region for the Stokes flow. 
\end{abstract}



\section{Introduction and statements of the main results}	

When two close-to-touching inclusions with high contrast material properties are present, the physical fields such as the stress may become arbitrarily large in the narrow region between them. Such field blow-up occurs in electro-statics and elasto-statics, and quantitative understanding of such a phenomenon is important in relation with the light confinement in the electro-static case, and with materials failure analysis in the elasto-static case. Lately, significant progress has been made in understanding the field enhancement. In the electro-static case, it is proved that the electric field, which is the gradient of the solution to the conductivity equation, blows up in the narrow region between two perfect conductors (where the conductivity is infinite) at the rate of $\Gd^{-1/2}$ \cite{AKL, Yun-SIAP-07} in two dimensions and of $|\Gd \log \Gd|^{-1}$ in three dimensions \cite{BLY-ARMA-09}, as the distance $\Gd$ between the two inclusions  tends to zero. The singular term of the stress concentration is also characterized in two dimensions \cite{ACKLY-ARMA-13}. This result has been extended to the elasticity in the context of the Lam\'e system of linear elasticity, showing that the blow-up rate of the stress in between two stiff inclusions (where the shear modulus is infinite) is $\Gd^{-1/2}$ in two dimensions \cite{BLL-ARMA-15, KY19}. References cited above are far from being complete. In fact, there is a long list of recent important achievements in this direction of research, for which we refer to the references in \cite{KY19,Milton-book-2001}.

In this paper, we consider the stress concentration in the two-dimensional steady Stokes system when two adjacent circular cylinders are present. Its quantitative analysis is important in understanding hydrodynamic interactions in soft matter systems. This problem is particularly interesting in comparison to the case of linear elasticity. In the linear elasticity case, the divergence of the displacement vector field blows up in general as the distance between two inclusions tends to zero, as was proved in \cite{KY19}. However, the divergence of the velocity vector in Stokes flow is confined to be zero, namely, the flow is incompressible. Thus, it is not clear whether the stress blows up or not in the case of Stokes flow, and how large it is if it actually blows up. The stress in the Newtonian fluid including the Stokes flow consists of two components, the pressure and the shear gradient of the velocity field. We investigate the blow-up rate of each component when the distance between the two cylinders tends to zero.

More precisely, suppose that two circular cylinders, denoted by $D_1$ and $D_2$, of the same radius $R$ are immersed in Stokes flow and they are separated by a distance $\Gd>0$. Since $D_1$ and $D_2$ are (rigid) cylinders, the boundary values of the steady flow on $\p D_1$ and $\p D_2$ are given as a linear combination of three vector fields representing rigid motions $\{\Bpsi_j\}_{j=1}^3$, which are defined as
\beq\label{Bpsi}
\Bpsi_1=\begin{bmatrix}
1 \\ 0
\end{bmatrix}, \quad
\Bpsi_2=\begin{bmatrix}
0 \\ 1
\end{bmatrix}, \quad
\Bpsi_3=\begin{bmatrix}
y \\ -x
\end{bmatrix}.
\eeq
Thus, we consider the following Stokes system in the exterior domain $D^e := \Rbb^2\setminus \overline{D_1\cup D_2}$:
\beq\label{stokes}
\begin{cases}
\mu \GD \Bu = \nabla p \quad &\mbox{in }D^e,
\\
\nabla \cdot \Bu = 0 \quad &\mbox{in }D^e,
\\[0.3em]
\ds \Bu = \sum_{j=1}^3 c_{ij} \Bpsi_j & \mbox{on } \p D_i, \  i=1,2,
\\[0.3em]
(\Bu-\BU, p-P) \in \Mcal_0,
\end{cases}
\eeq
where $\Gm$ represents the constant viscosity of the fluid, $c_{ij}$ are constants to be determined from the equilibrium conditions (see \eqnref{equili} below), $(\BU,P)$ is a given background solution to the homogeneous Stokes system in $\Rbb^2$, namely,
\beq
\mu \GD \BU = \nabla P \quad \mbox{in } \Rbb^2,
\eeq
and the class $\Mcal_0$ is characterized by decay conditions at $\infty$. The precise definition of $\Mcal_0$ is given later in Subsection \ref{subsec:diri}. Here we just mention that the problem \eqnref{stokes} admits a unique solution.

Throughout this paper, we assume that both the gradient $\nabla \BU$ of the background velocity field and the pressure $P$ are constant functions. Since the  pressure is determined up to a constant, we assume that $P=0$ and
\beq\label{fieldU}
\BU(x,y) = \begin{bmatrix} a & c \\ d & -a \end{bmatrix} \begin{bmatrix} x \\ y \end{bmatrix} \quad (a^2 + (c+d)^2 \neq 0)
\eeq
for some constants $c$ and $d$. The fields in \eqnref{fieldU} are the only divergence-free fields in the case where $\nabla \BU$ is constant. The condition $a^2 + (c+d)^2 \neq 0$ is imposed from the fact that otherwise
$$
\BU(x,y) = c\begin{bmatrix} y \\ -x \end{bmatrix} ,
$$
and hence $\BU$ with a constant $p$ is the solution to the problem \eqnref{stokes}, and its gradient does not blow up. If we write $\BU$ as
\beq
\BU = a \BUst + \frac{c+d}{2} \BUsh + \frac{c-d}{2} \begin{bmatrix} y \\ -x \end{bmatrix} := a \begin{bmatrix} x \\ -y \end{bmatrix} + \frac{c+d}{2} \begin{bmatrix} y \\ x \end{bmatrix} + \frac{c-d}{2} \begin{bmatrix} y \\ -x \end{bmatrix},
\eeq
and denote respectively by $(\Bust, \pst)$ and $(\Bush, \psh)$ the solutions to \eqnref{stokes} when $\BU=\BUst$ and $\BU=\BUsh$, then the solution $(\Bu, p)$ is given by
\beq\label{Buexpre}
\Bu = a \Bust + \frac{c+d}{2} \Bush + \frac{c-d}{2} \begin{bmatrix} y \\ -x \end{bmatrix}
\eeq
and
\beq\label{pexpre}
p = a \pst + \frac{c+d}{2} \psh + \frac{c-d}{2} \mbox{const.}
\eeq
The singular behavior of the stress comes solely from those corresponding to $(\Bust,\pst)$ and $(\Bush, \psh)$. The flows $\BU_{\rm ex}$ and $\BU_{\rm sh}$ are called the extensional flow and the shear flow, respectively, which explains the subscripts ex and sh in our notation.

For the solution $(\Bu, p)$ to the Stokes system, the strain tensor, denoted by $\Ecal[\Bu]$, is given by
\beq
\Ecal[\Bu] = \frac{1}{2}(\nabla \Bu + \nabla\Bu^T) ,
\eeq
where the superscript $T$ denotes the transpose, and the corresponding stress tensor is given by
\beq\label{Gsdef}
\Gs[\Bu,p] = - p I + 2\mu \Ecal[\Bu],
\eeq
where $I$ is the identity matrix. The constants $c_{ij}$ appearing in \eqnref{stokes} are determined by the boundary integral conditions
\beq\label{equili}
\int_{\p D_i} \Bpsi_j \cdot \Gs[\Bu,p] \nu \, dl=0, \quad  i=1,2, \ j=1,2,3.
\eeq
Here, $\nu$ denotes the unit normal on the boundary $\p D_i$ and $dl$ is the line element.
Physically, these integral conditions imply that each rigid inclusion is in equilibrium, namely, the net translational and rotational stress on each boundary is zero (see, e.g., \cite{Berlyand-SIMA-06}).

The following is the main result of this paper. It shows that the stress always blows up. There and in what follows, $A \lesssim B$ means that there is a constant $C$ independent of $\Gd$ such that $A \le C B$, and $A \approx B$ means that both $A \lesssim B$ and $B \lesssim A$ hold. The supremum norm on $D^e$ is denoted by $\| \cdot \|_\infty$.

\begin{theorem}\label{thm:main}
Let $D_1$ and $D_2$ be disks of the same radii and let $(\Bu,p)$ be the the unique solution to \eqnref{stokes} when $U$ is of the form \eqnref{fieldU} and $P=0$. Then,
\beq
\| \Gs[\Bu,p] \|_\infty \approx \Gd^{-1/2}.
\eeq
\end{theorem}

In fact, we can separate our problem into the cases where the pressure or shear stress blows up as the following two theorems show, of which the main theorem is an immediate consequence. To present these results clearly, we assume for convenience that the centers of $D_1$ and $D_2$ are, respectively, given by
\beq\label{config}
c_1= (-R-\Gd/2,0) \quad\mbox{and}\quad c_2=(R+\Gd/2,0)
\eeq
after applying rotation and translation if necessary, where $R$ is the common radius of the disks and $\Gd$ is the distance between them.  To describe the two-dimensional Stokes flow, we construct a pair of stream functions using the bipolar coordinates, and then use the stream function formulation to construct special solutions $(\Bh_j,p_j)$, $j=1,2$, to the Stokes system (see Section \ref{sec:singular} for precise definitions of $(\Bh_j,p_j)$). It turns out that these special solutions, called singular functions, capture precisely the singular behavior of $\Gs[\Bust,\pst]$ and $\Gs[\Bush,\psh]$. As a result, we are able to characterize the blow-up of the pressure and the shear stress for the different configurations of the background velocity field $\BU$: when $\BU=\BUst$, the pressure blows up at the rate of $\Gd^{-1/2}$ while the shear stress is bounded; when $\BU=\BUsh$, the other way around.

The precise statements of the results are presented in the following theorems. Here and afterwards, $\Pi_\Gd$ denotes the narrow region  between the two cylinders defined by
\beq\label{narrow}
\Pi_\Gd := ([-R-\Gd/2, R+\Gd/2] \times [-\sqrt{\Gd}, \sqrt{\Gd}]) \cap D^e.
\eeq

\begin{theorem}\label{thm:main1}
Suppose that $D_1$ and $D_2$ are arranged so that \eqnref{config} holds and that $\BU=\BUst$ and $P=0$. It holds that
\beq
\| \Ecal[\Bust] \|_\infty \lesssim 1 \quad\mbox{and}\quad \| \pst \|_\infty \approx \Gd^{-1/2}
\eeq
as $\Gd \to 0$. In the narrow region $\Pi_\Gd$,
\beq\label{stressst}
\Gs[\Bust,\pst](x,y) = 2\Gm \sqrt{R} \Gd^{-1/2} \frac{(y^2+3R\Gd)(y^2-R\Gd)}{(y^2+R\Gd)^2} I + O(1).
\eeq
\end{theorem}

\begin{theorem}\label{thm:main2}
Suppose that $D_1$ and $D_2$ are arranged so that \eqnref{config} holds and that $\BU=\BUsh$ and $P=0$. It holds that
\beq
\|\Ecal[\Bush]\|_\infty \approx \Gd^{-1/2} \quad\mbox{and}\quad
\| \psh \|_\infty \lesssim 1
\eeq
as $\Gd \rightarrow 0$. In the narrow region $\Pi_\Gd$,
\beq\label{stresssh}
\Gs[\Bush,\psh](x,y)=    2\mu  \sqrt{\frac{R}{\Gd}} \frac{R\Gd}{y^2+R\Gd} \begin{bmatrix}0 & 1 \\ 1 & 0 \end{bmatrix} + O(1).
\eeq
\end{theorem}

Let $(\Bu,p)$ be the solution to \eqnref{stokes}. According to \eqnref{Buexpre} and \eqnref{pexpre},
\beq
\Gs[\Bu,p] = a \Gs[\Bust,\pst] + \frac{c+d}{2} \Gs[\Bush,\psh] +O(1).
\eeq
Thus, Theorem \ref{thm:main} is an immediate consequence of \eqnref{stressst} and \eqnref{stresssh}.

What is actually shown in this paper is that if the background velocity field $\BU$ is of the form \eqnref{fieldU} and  $P=0$, then the solution $(\Bu, p)$ is of the following form:
\beq\label{Bpcomp}
(\Bu,p) = a \frac{2}{\sqrt{R}} \Gd^{3/2} (\Bh_1,p_1) + \frac{c+d}{2} \sqrt{R\Gd} (\Bh_2, p_2) + (\Bu_0, p_0),
\eeq
where $(\Bu_0, p_0)$ is a solution to the Stokes problem whose stress tensor is bounded. See the end of section \ref{sec:1.2} for a brief proof of this fact. Since the singular functions $(\Bh_j,p_j)$ ($j=1,2$) are given explicitly, the decomposition formula \eqnref{Bpcomp} may cast light on the challenging problem of computing the Stokes flow in presence of closely located rigid cylinders.

Some historical remarks on the study of the Stokes flow in presence of two circular cylinders are in order. Jeffrey developed in \cite{Jeffrey-PTRS-1921} a separable solution method based on bipolar coordinates and then analyzed in \cite{Jeffrey-PRSA-1922} the flow generated by two rotating circular cylinders. Several other authors independently developed similar methods \cite{Wannier-QAM-1950, BE-PF-1965}.  Jeffrey's method has been applied to various problems of Stokes flow \cite{Schubert-JFM-1967, Wakiya-JPSJ-1975, Wakiya-JPSJ-1975-II, JO-QJMAM-1981, Smith-M-1991, Watson-M-1995, Crowdy-IJNLM-2011, IC-JFM-2017}. In particular, Raasch derived the exact analytic solution for two circular cylinders under the equilibrium condition, which represents suspended particles in a viscous fluid \cite{Raasch-PhD-1961} (see also \cite{Raasch-ZAMM-1961, DRM-CJCE-1967}). However, due to the high complexity of the solution, it is difficult to analyze the singular behavior of the solution when the cylinders are close-to-touching. In this work, this difficulty is successfully overcome by introducing the singular functions.

Other than the method of bipolar coordinates, a formal asymptotic technique called the lubrication theory was also developed for the viscous flow in the narrow region \cite{FA-CES-1967,Graham-ASR-1981,NK-JFM-1984}. Berlyand {\it et al} \cite{Berlyand-SIMA-06} constructed a refined lubrication approximation and then derived an asymptotic formula for the effective viscosity of concentrated suspensions. We mention that the approximation \eqref{Bpcomp} is different from the lubrication one in two respects. Firstly, it provides a rigorous pointwise approximation of the solution in the narrow region. Secondly, its singular parts satisfy the Stokes system at the exact level, which is a key to the development of an accurate numerical scheme.

The organization of the paper is as follows. In the next section, we introduce the bipolar coordinates and review the stream function formulation for the Stokes system. Section \ref{sec:singular} is to construct singular functions which are the building blocks in describing the singular behavior of the solution to the Stokes system \eqnref{stokes} as the separating distance between $D_1$ and $D_2$ tends to zero. Sections \ref{sec:1.1} and \ref{sec:1.2} are to prove Theorems \ref{thm:main1} and \ref{thm:main2}. Section \ref{sec:noslip1} and \ref{sec:noslip2} are to prove that stress does not blow up if the no-slip boundary condition is prescribed on the boundary of the circular inclusions. Appendices are to prove some auxiliary lemmas. The paper ends with a short discussion.

\section{Preliminaries}

\subsection{Bipolar coordinates}\label{subsec:bi}

Given a positive constant $a$, the bipolar coordinates $(\Gz,\Gt)$ are defined by
\beq
x + i y = a \frac{e^{\Gz - i\Gt}+1}{e^{\Gz - i \Gt}-1},
\eeq
so that
\beq\label{eq:bipolar_x_y}
x = a \frac{\sinh \Gz}{\cosh \Gz - \cos\Gt}, \quad  y= a \frac{\sin \Gt}{\cosh \Gz - \cos\Gt},
\eeq
or equivalently,
\beq\label{cartbipolar}
\Gz = \log \frac{\sqrt{(x+a)^2+y^2}}{\sqrt{(x-a)^2+y^2}}, \quad \Gt= \arg(x-a,y)-\arg(x+a,y).
\eeq
The coordinate curve $\{\Gz = c\}$ represents a circle of radius $a/|\sinh c|$ centered at the point $(a/\tanh c,0)$. Similarly, the curve  $\{\Gt = c\}$ represents a circle of radius $a/|\sin c|$ centered at $(a/\tan c,0)$.
The point of infinity corresponds to $(\Gz,\Gt)=(0,0)$. See, e.g., \cite{Jeffrey-PTRS-1921} for bipolar coordinates in relation with the Stokes system.

The geometry of two disks (the cross sections of the two circular cylinders) can be described efficiently in terms of bipolar coordinates. Let
\beq\label{adef}
a := \sqrt{\Gd \left( R + \frac{\Gd}{4} \right)}.
\eeq
Then the boundary $\p D_i$ of the cylinder $D_i$ can be parameterized by a $\Gz$-coordinate curve as follows:
\beq
\p D_1 = \{ \Gz = - s\}, \quad \p D_2 = \{ \Gz = + s\},
\eeq
where
\beq\label{s}
s = \sinh^{-1}(a/R).
\eeq
We note that
\beq\label{sGd}
s = \sqrt{\frac{\Gd}{R}} + O(\Gd^{3/2}) \quad\mbox{as } \Gd \to 0.
\eeq
The exterior domain $D^e$ of $D_1 \cup D_2$ is characterized in bipolar coordinates $(\Gz, \Gt)$ by the rectangle
\beq
D^e= \{ (\Gz, \Gt) \in (-s,s) \times [0, 2\pi) \}.
\eeq
In particular, $\{ (\Gz, \pi) , |\Gz| <s \}$ is the line segment connecting the two points $(-\Gd/2, 0)$ and $(\Gd/2,0)$.

Let $\{\Be_x,\Be_y \}$ be the standard unit basis vectors in $\Rbb^2$ and let $\{\Be_{\Gz},\Be_\Gt \}$ be the unit basis vectors in the bipolar coordinates, namely,
$$
\Be_\Gz = \frac{\nabla \Gz}{|\nabla \Gz|}, \quad
\Be_\Gt = \frac{\nabla \Gt}{|\nabla \Gt|}.
$$
Let $[\Be_\Gz, \Be_\Gt]$ denote the $2\times 2$ matrix whose columns are $\Be_\Gz$ and $\Be_\Gt$. Then one can easily see from \eqnref{cartbipolar} that
\beq\label{Xione}
\Xi := [\Be_\Gz, \Be_\Gt] = \begin{bmatrix} \Ga(\Gz,\Gt) & -\Gb(\Gz,\Gt) \\ - \Gb(\Gz,\Gt) & - \Ga(\Gz,\Gt) \end{bmatrix},
\eeq
where
\beq\label{pqdef}
\Ga(\Gz,\Gt) := \frac{1-\cosh\Gz\cos\Gt  }{\cosh\Gz - \cos\Gt} , \quad \Gb(\Gz,\Gt) := \frac{\sinh\Gz \sin\Gt}{\cosh\Gz-\cos\Gt}.
\eeq
Since $\Ga^2+\Gb^2=1$, we have
\beq\label{Xisquare}
\Xi^2=I.
\eeq
This means that $\Xi$ is the transition transformation in the sense that
\beq\label{exey}
[\Be_x, \Be_y] = \Xi [\Be_\Gz, \Be_\Gt].
\eeq

Define the scaling function
\beq\label{h_def}
h(\Gz,\Gt) : = \frac{\cosh\Gz-\cos\Gt}{a}.
\eeq
Then, for any scalar function $f$, its gradient $\nabla f$ can be expressed as
\beq\label{gradrel}
\nabla f = h(\Gz,\Gt) [\p_\Gz f \,\Be_\Gz + \p_\Gt f \,\Be_\Gt]
\eeq
(see, e.g., \cite{Smythe-1968}). Here and throughout this paper, $\p_\Gz$ and $\p_\Gt$ denote the partial derivatives with respect to the $\Gz$ and $\Gt$ variables, respectively. Moreover, the line element, denoted by $dl$, on $\p D_2$ is given by
\beq\label{lineele}
dl = h(s,\Gt)^{-1} d\Gt.
\eeq
One can easily check that, for $i=1,2$,
\beq\label{eq:nor_deri_bipolar}
\p_\nu f\big|_{\p D_i} = (-1)^{i+1} h(\Gz,\Gt) \p_\Gz f\big|_{\Gz=(-1)^i s} \ ,
\eeq
and
\beq\label{eq:tan_deri_bipolar}
\p_T f\big|_{\p D_i} = (-1)^i h(\Gz,\Gt) \p_\Gt f\big|_{\Gz=(-1)^i s} \ ,
\eeq
where $\p_\nu$ and $\p_T$ denote the normal and tangential derivatives, respectively.

Using \eqref{eq:bipolar_x_y} and \eqnref{adef}, one can see that
$$
\frac{\cos \Gt}{\cosh\Gz - \cos \Gt} = \frac{1}{2a^2}(x^2+y^2) -\frac{1}{2} = \frac{1}{2R\Gd}(x^2+y^2) -\frac{1}{2} + O(\Gd).
$$
If $(x,y)$ lies in the narrow region $\Pi_\Gd$ defined in \eqnref{narrow}, then $|x| \lesssim \Gd$, and hence
$$
\frac{1}{2R\Gd}(x^2+y^2) -\frac{1}{2} = \frac{y^2}{2R\Gd} -\frac{1}{2} +O(\Gd).
$$
Moreover, if $(\Gz,\Gt)$ lies in $\Pi_\Gd$, then there is a positive constant $C< \pi$ such that $|\Gt-\pi| < C$. Since $|\Gz| < s \approx \sqrt{\Gd}$, we have
$$
\cosh\Gz - \cos \Gt = 1 - \cos \Gt + O(\Gz^2)= 1 - \cos \Gt + O(\Gd).
$$
Thus we have
$$
\frac{\cos \Gt}{1 - \cos \Gt} = \frac{y^2}{2R\Gd} -\frac{1}{2} + O(\Gd),
$$
or equivalently,
\beq\label{cosnarrow}
\cos\Gt= \frac{y^2-R\Gd}{y^2+R\Gd} + O(\Gd)
\eeq
in the region $\Pi_\Gd$.
One can also easily see from \eqnref{pqdef} that in $\Pi_\Gd$
$$
\Ga(\Gz,\Gt) = 1+ O(\Gd), \quad \Gb(\Gz,\Gt)=O(\sqrt{\Gd}),
$$
and hence
\beq\label{Xinarrow}
\Xi =\begin{bmatrix} 1 & 0 \\ 0 & - 1 \end{bmatrix} + O(\sqrt{\Gd}).
\eeq

Using \eqref{eq:bipolar_x_y}, one can see
$$
|\Bx|^2 = x^2+y^2 = \frac{\cosh\Gz+\cos\Gt}{\cosh\Gz-\cos\Gt}.
$$
Since the following relation holds for large enough $|\Bx|$ (or small enough $\Gz$ and $\Gt$)
$$
|\Bx|^{-2}=\frac{\cosh\Gz-\cos\Gt}{\cosh\Gz+\cos\Gt} = \frac{\frac{\Gz^2}{2} + \frac{\Gt^2}{2}}{2+O(\Gz^2+\Gt^2)},
$$
we obtain
\beq\label{eq:largex_smallGzGt}
\frac{1}{8}({\Gz^2+\Gt^2}) \leq |\Bx|^{-2} \leq \frac{1}{2}(\Gz^2+\Gt^2).
\eeq

\subsection{The stream function}

Here we review the stream function formulation in the two-dimensional incompressible flow and collect some useful formulas.

It is well known that any solution $(\Bu,p)$ to the Stokes system, $\mu \GD \Bu = \nabla p$ and $\nabla \cdot \Bu = 0$, can be written using a scalar function $\Psi$ satisfying the biharmonic equation $\GD^2 \Psi = 0$. The function $\Psi$ is called the stream function. Once the function $\Psi$ is known, the velocity field $\Bu=(u_x,u_y)^T$ can be determined from the relations
\beq\label{uxPsi}
u_x = \p_y \Psi, \quad u_y = -\p_x \Psi,
\eeq
and the pressure $p$ is a harmonic conjugate of $\mu \GD \Psi$ (see, e.g., \cite{Batchelor-1967}).

Let us write the stream function formulation in terms of bipolar coordinates. It is also known (see, e.g, \cite{Wakiya-JPSJ-1975, Wakiya-JPSJ-1975-II}) that the Laplacian in Cartesian coordinates is related to bipolar coordinates via
\beq\label{eq:Laplacian_Psi}
\GD_{x,y} \Psi = \frac{1}{a}\left( (\cosh\Gz-\cos\Gt) \GD_{\Gz,\Gt}  + (\cosh\Gz + \cos\Gt) - 2(\sinh\Gz \p_\Gz + \sin\Gt \p_\Gt) \right)(h\Psi),
\eeq
where $h$ is the function defined in \eqnref{h_def}. Using this formula, the biharmonic equation $\GD^2\Psi=0$ can be rewritten as
\beq\label{eq:biharmonic_bipolar}
\left( \p_\Gz^4 + 2 \p_\Gz^2 \p_\Gt^2 +  \p_\Gt^4 -2 \p_\Gz^2 + 2 \p_\Gt^2 + 1 \right) (h\Psi) =0,
\eeq
and the general solution to the above equation takes the following form:
\begin{align}
& (h\Psi)(\Gz,\Gt) \nonumber \\
& = K(\cosh\Gz -\cos\Gt)\ln (2\cosh \Gz - 2\cos \Gt)
+ a_0 \cosh \Gz + b_0 \Gz \cosh \Gz + c_0 \sinh \Gz + d_0 \Gz \sinh \Gz \nonumber\\
& + (a_1 \cosh 2\Gz + b_1 + c_1 \sinh 2\Gz + d_1 \Gz)\cos\Gt \nonumber
 + (\widetilde{a}_1 \cosh 2\Gz + \widetilde{b}_1 + \widetilde{c}_1 \sinh 2\Gz + \widetilde{d}_1 \Gz)\sin\Gt \nonumber
\\
&
+\sum_{n=2}^\infty \Big(a_n \cosh (n+1) \Gz + b_n \cosh (n-1)\Gz + c_n \sinh (n+1)\Gz+ d_n \sinh(n-1)\Gz \Big)\cos n\Gt\nonumber
\\
&
+\sum_{n=2}^\infty \Big(\widetilde{a}_n \cosh (n+1) \Gz + \widetilde{b}_n \cosh (n-1)\Gz +\widetilde{c}_n \sinh (n+1)\Gz+ \widetilde{d}_n \sinh(n-1)\Gz \Big)\sin n\Gt.
\label{eq:general_biharmonic_solution}
\end{align}

Using \eqnref{gradrel} and \eqnref{uxPsi}, one can see that the components of the velocity $\Bu = u_\Gz \Be_\Gz + u_\Gt \Be_\Gt$ are given as follows:
\begin{align}
u_\Gz &=  - h \p_\Gt \Psi =\left(- \p_\Gt + \frac{\sin\Gt}{\cosh\Gz-\cos\Gt}\right)(h\Psi), \label{eq_velo1} \\
u_\Gt &= +h \p_\Gz \Psi = \left( \p_\Gz - \frac{\sinh\Gz}{\cosh\Gz-\cos\Gt}\right)(h\Psi), \label{eq_velo2}
\end{align}
and the pressure $p$ satisfies the relations
\beq\label{eq:pressure_bipolar}
\p_\Gz p = - \mu \p_\Gt \GD \Psi, \quad
\p_\Gt p = \mu \p_\Gz  \GD \Psi.
\eeq

The entries of the strain tensor $\Ecal[\Bu]$ when represented in terms of the basis $\{\Be_\Gz, \Be_\Gt\}$ are given by
\begin{align}
\Ecal_{\Gz\Gz}&= - h \p_\Gz \left( h \p_\Gt \Psi \right) - h \p_\Gz \Psi \p_\Gt h,
\label{eq:strain_bipolar1}
\\
\Ecal_{\Gt\Gt}&= + h \p_\Gt \left( h \p_\Gz \Psi  \right) + h \p_\Gt \Psi \p_\Gz h,
\label{eq:strain_bipolar2}
\\
\Ecal_{\Gz\Gt}&= \frac{1}{2} \left( \p_\Gz  \left( {h^2} \p_\Gz \Psi  \right) - \p_\Gt \left({h^2} \p_\Gt \Psi  \right)\right).
\label{eq:strain_bipolar3}
\end{align}
Therefore, the following relation holds:
\beq\label{Ecalrel}
\Ecal[\Bu]= \Xi \begin{bmatrix} \Ecal_{\Gz\Gz} & \Ecal_{\Gz\Gt} \\ \Ecal_{\Gz\Gt} & \Ecal_{\Gt\Gt} \end{bmatrix} \Xi,
\eeq
where $\Xi$ is the matrix given in \eqnref{Xione}. The entries of the stress tensor in bipolar coordinates are given by
\beq\label{eq:sigma_formula_bipolar}
\Gs _{\Gz\Gz} = - p + 2\mu \Ecal_{\Gz\Gz},
\quad
\Gs_{\Gt\Gt} = - p + 2\mu \Ecal_{\Gz\Gt},
\quad
\Gs_{\Gz\Gt} = 2\mu \Ecal_{\Gz\Gt}.
\eeq
Similarly, we have the following relation for the stress tensor:
\beq\label{Gsrel}
\Gs[\Bu,p]= \Xi \begin{bmatrix} \Gs_{\Gz\Gz} & \Gs_{\Gz\Gt} \\ \Gs_{\Gz\Gt} & \Gs_{\Gt\Gt} \end{bmatrix} \Xi.
\eeq
Since each component of $\Xi$ is bounded, it follows from \eqnref{Xisquare}, \eqnref{Ecalrel} and \eqnref{Gsrel} that
\beq\label{Ecalnorm}
\| \Ecal[\Bu] \|_{L^\infty(K)} \approx \| \Ecal_{\Gz\Gz} \|_{L^\infty(K)} +  \| \Ecal_{\Gt\Gt} \|_{L^\infty(K)} + \| \Ecal_{\Gz\Gt} \|_{L^\infty(K)},
\eeq
and
\beq\label{Gsnorm}
\|\Gs[\Bu,p] \|_{L^\infty(K)} \approx \| \Gs_{\Gz\Gz} \|_{L^\infty(K)} +  \| \Gs_{\Gt\Gt} \|_{L^\infty(K)} + \| \Gs_{\Gz\Gt} \|_{L^\infty(K)}
\eeq
for any subset $K$ of $D^e$.

Using integrations by parts on the exterior domain $D^e$, we have for any solutions $(\Bu,p),(\Bv,q)$ to the Stokes system such that $\Bu(\Bx),\Bv(\Bx)\rightarrow 0$ as $|\Bx|\rightarrow\infty$ that
\begin{align}
\int_{ \p D^e}  \Bu\cdot\Gs[\Bv,q]\nu &= -\int_{ D^e}  \Ecal[\Bu]:\Gs[\Bv,q]\nonumber
\\
&=\int_{ D^e}  (\nabla \cdot \Bu)q -2\mu\int_{ D^e}  \Ecal[\Bu]:\Ecal[\Bv]\nonumber
\\
&= -2\mu\int_{ D^e}  \Ecal[\Bu]:\Ecal[\Bv].\label{eq:int_parts_formula}
\end{align}
This implies in particular that the following Green's theorem holds:
\beq\label{eq:int_parts_formula2}
    \int_{ \p D^e}  \Bu\cdot\Gs[\Bv,q]\nu = \int_{ \p D^e}  \Bv\cdot\Gs[\Bu,p]\nu.
\eeq

\subsection{An exterior Dirichlet problem}\label{subsec:diri}

Let $\BGG(\Bx)=(\GG_{ij}(\Bx))_{i,j=1,2}$ be
$$
\GG_{ij}(\Bx) = -\frac{1}{4\pi\mu}(\delta_{ij}\log|\Bx| + \frac{x_i x_j}{|\Bx|^2}), \quad \Bx\in\Rbb^2\setminus\{(0,0)\},
$$
and define $\Bp=(p_j)_{j=1,2}$ by
$$
\Bp = -\frac{1}{2\pi} \frac{\Bx}{|\Bx|^2}, \quad \Bx\in\Rbb^2\setminus\{(0,0)\}.
$$
Then, $(\BGG,\Bp)$ is the fundamental solution to the Stokes system, namely,
$$
\mu \GD \BGG - \nabla \Bp = \delta(\Bx)\mathbf{I}.
$$
Let $\GG_\GD$ be the fundamental solution to the Laplacian given by
$$
\GG_\GD (\Bx) = \frac{1}{2\pi}\log|\Bx|.
$$

The existence and uniqueness of the exterior Dirichlet problem, proved in \cite[Theorem 9.15]{Mitrea-book-2012}, is as follows.

\begin{theorem}\label{thm:ext_diri}
Assume that $\GO$ is a bounded Lipschitz domain. Then the exterior Dirichlet problem
\beq
\begin{cases}
\mu \GD \Bu = \nabla p &\quad \mathrm{in } \ \Rbb^2\setminus \overline{\GO},
\\
\nabla \cdot \Bu = 0 &\quad \mathrm{in }\ \Rbb^2\setminus \overline{\GO},
\\[0.3em]
\Bu = \Bg &\quad \mathrm{on }\ {\p \GO},
\end{cases}
\eeq
with the decaying conditions
\begin{align*}
  \begin{cases}
 \Bu(\Bx) = \mathbf{\GG}(\Bx)\BA +\BC+O(|\Bx|^{-1}), \\
\p_j \Bu(\Bx) = \p_j \mathbf{\GG}(\Bx)\BA +O(|\Bx|^{-2}), \\
p(\Bx) = \nabla \GG_\GD \cdot \BA + O(|\Bx|^{-2})
\end{cases}
\end{align*}
as $|\Bx|\rightarrow\infty$ for some constant $\BC\in\Rbb^2$, has a solution, which is unique modulo adding functions to the
pressure term which are locally constant in $\Rbb^2$.
Here, $\BA \in \Rbb^2$ is a priori given constant.

\end{theorem}

We shall consider the exterior Dirichlet problem with $\BA=0$.
Let $\Mcal$ be the set of all pairs of functions $(\Bu,p)$ satisfying
\beq\label{Mcal}
  \begin{cases}
 \Bu(\Bx) = \BC+ O(|\Bx|^{-1}),\\
\nabla \Bu(\Bx) =  O(|\Bx|^{-2}), \\
p(\Bx) = O(|\Bx|^{-2})
\end{cases}
\eeq
as $|\Bx|\rightarrow\infty$ for some constant $\BC\in\Rbb^2$. We denote by $\Mcal_0$ the set of all pairs of functions $(\Bu,p)$ satisfying the decay conditions \eqref{Mcal} with $\BC=0$.

\section{The singular functions for the Stokes system}\label{sec:singular}

In what follows, we construct the singular functions $(\Bh_j,p_j)$, $j=1,2$, which is the unique solution to the following problem:
\beq\label{Bhj}
\begin{cases}
\mu \GD \Bh_j = \nabla p_j &\quad \mbox{in }D^e,
\\
\nabla \cdot \Bh_j = 0 &\quad \mbox{in } D^e,
\\[0.3em]
\Bh_j = \frac{(-1)^i}{2}\Bpsi_j&\quad {\p D_i}, \ i=1,2,
\\[0.3em]
(\Bh_j,p_j) \in \Mcal.
\end{cases}
\eeq
We then provide quantitative estimates of the blow-up of these functions in the subsequent propositions. We call the solutions $(\Bh_j,p_j)$ the singular functions since they are the building blocks in describing the singular behavior, i.e., the stress tensor blow-up, of the solution to \eqnref{stokes}. In fact, we will see that the solution to \eqnref{stokes} can be expressed as a linear combination of singular functions (modulo a regular function) and the nature of the stress tensor blow-up is characterized by $(\Bh_1,p_1)$ and $(\Bh_2,p_2)$.

\begin{prop}\label{lem:hone}
Define two constants $A_1$ and $B_1$ by
\beq\label{AB}
A_1 := \frac{1}{2s - \tanh 2s}, \quad
B_1 := -\frac{1}{2\cosh 2s}A_1.
\eeq

\begin{itemize}
\item[{\rm (i)}] The stream function  $\Psi_1$ associated with the singular functions $(\Bh_1,p_1)$ is given by
\beq\label{eq:stream_singular1}
\Psi_1(\Gz,\Gt) = \frac{1}{h(\Gz,\Gt)}( A_1 \Gz + B_1 \sinh 2\Gz) \sin\Gt.
\eeq

\item[{\rm (ii)}] The components of the velocity $\Bh_1 = h_{1\Gz}\Be_\Gz + h_{1\Gt} \Be_\Gt$ are given by
\begin{align}
h_{1\Gz} &= (A_1\Gz + B_1\sinh 2\Gz) \frac{1-\cosh\Gz\cos\Gt}{\cosh\Gz-\cos\Gt}, \label{honeone}
\\
h_{1\Gt} &= \sin\Gt \left( A_1+2B_1 \cosh 2\Gz - \frac{\sinh\Gz (A_1\Gz + B_1\sinh 2\Gz)}{\cosh\Gz-\cos\Gt} \right). \label{honetwo}
\end{align}

\item[{\rm (iii)}]
The pressure $p_1$  is given by
\beq\label{eq:pressure1}
p_1 = \frac{2\mu}{a}( (A_1-2B_1)\cosh\Gz \cos\Gt + B_1 \cosh 2\Gz \cos 2\Gt) - \frac{2\mu}{a}(A_1-B_1).
\eeq

\end{itemize}
\end{prop}

\proof
The formulas \eqnref{eq:stream_singular1}-\eqnref{honetwo} are derived in the following way. We use the expansion \eqref{eq:general_biharmonic_solution} for the general solution to the Stokes system, and then determine its unknown constant coefficients by matching the boundary conditions on $\p D^e$, given by $\{\Gz=\pm s\}$, and using formulas \eqref{exey}, \eqref{eq_velo1}, and \eqref{eq_velo2}. Let us show that the boundary conditions are fulfilled. If $\Gz = \pm s $, we have
\begin{align*}
h_{1\Gz} |_{\Gz=\pm s} &= \pm \frac{1-\cosh s\cos\Gt  }{2(\cosh s - \cos\Gt)}= \pm \frac{1}{2}\Ga(s,\Gt), \\
h_{1\Gt} |_{\Gz=\pm s} &= \mp \frac{ \sinh s \sin\Gt}{2(\cosh s -\cos\Gt)} = \mp \frac{1}{2}\Gb(s,\Gt) .
\end{align*}
One can see from the relation \eqref{exey} that the boundary conditions on $\p  D_1\cup \p D_2$ are satisfied.

The formula \eqnref{eq:pressure1} follows from \eqref{eq:pressure_bipolar} and \eqnref{eq:stream_singular1}. In fact, applying \eqref{eq:Laplacian_Psi} to $\Psi_1$ given in \eqnref{eq:stream_singular1}, we see that
$$
\mu\GD \Psi_1 = \frac{-2\mu}{a}( (A_1-2B_1)\sinh\Gz \sin\Gt + B_1 \sinh 2\Gz \sin 2\Gt).
$$
The harmonic conjugate of this function vanishing at $(\Gz,\Gt)=(0,0)$ is nothing but the one given in \eqnref{eq:pressure1}.

We now prove that $(\Bh_1, p_1)$ belongs to $\Mcal$. We first prove that $\Bh_1(\Bx) = O(|\Bx|^{-1})$ as $|\Bx|\rightarrow\infty$, which amounts to proving
\beq\label{eq:h1_far_claim}
\Bh(\Gz,\Gt) = O(|\Gz|+|\Gt|), \quad (\Gz,\Gt)\rightarrow(0,0),
\eeq
thanks to \eqref{eq:largex_smallGzGt}.
We have from \eqref{honeone} and \eqref{honetwo} that
\begin{align}
|h_{1\Gz}|&\leq C (|\Gz| 	+|\Gz| \frac{|\Gz|^2+|\Gt|^2}{|\Gz|^2+|\Gt|^2}) \leq C |\Gz|,\nonumber
\\
|h_{1\Gt}|&\leq C |\Gt|\Big(1 	+|\Gz| \frac{|\Gz|}{|\Gz|^2+|\Gt|^2}\Big) \leq C |\Gt|.\label{eq:h1GzGt_far_estim}
\end{align}
Here and throughout this proof, the constant $C$ may depend on $s$, but is independent of $(\Gz,\Gt)$.
This proves \eqref{eq:h1_far_claim}.

Similarly, one can show that
\begin{align}
	|\p_\Gz h_{1\Gz}|\leq C, \quad |\p_\Gz h_{1\Gz}|\leq C, \quad |\p_\Gz h_{1\Gz}|\leq C, \quad |\p_\Gz h_{1\Gz}|\leq C.
	\label{eq:h1GzGt_grad_far_estim}
\end{align}
 Since $\Bh_1 = h_{1\Gz}\Be_\Gz + h_{1\Gt} \Be_\Gt$, we have
\begin{align}
|\nabla \Bh_1| \leq C (|\nabla h_{1\Gz}| + |h_{1\Gz} \nabla \Be_\Gz| + |\nabla h_{1\Gt}| + |h_{1\Gt} \nabla \Be_\Gt|).
\end{align}
It then follows from \eqref{eq:h1GzGt_far_estim} and the following lemma, whose proof will be given in Appendix \ref{appendixC}, that
\begin{align}
|\nabla \Bh_1| &\leq C (|\nabla h_{1\Gz}| + |\nabla h_{1\Gt}| + |\Gz|^2 +|\Gt|^2).
\end{align}

\begin{lemma}\label{lem:grad_eGz_eGt_estim}
It holds that
\beq
|\nabla \Be_\Gz| + |\nabla \Be_\Gt| \lesssim |\Gz| + |\Gt|.
\eeq
\end{lemma}

We then have from \eqref{eq:h1GzGt_grad_far_estim} that
\begin{align*}
|\nabla \Bh_1|
\leq  C (| h \p_\Gz h_{1\Gz} | +	| h \p_\Gt h_{1\Gz} |+| h \p_\Gz h_{1\Gt} |+| h \p_\Gt h_{1\Gt} |+|\Gz|^2+|\Gt|^2)
\leq C ( |h| + |\Gz|^2 + |\Gt|^2) .
\end{align*}
One can see from the definition of the function $h$ that
$$
|h(\Gz,\Gt)|  \leq C (|\Gz|^2 + |\Gt|^2),
$$
and hence
$$
|\nabla \Bh_1| \leq C (|\Gz|^2 + |\Gt|^2),
$$
or equivalently, $\nabla \Bh_1(\Bx) = O(|\Bx|^{-2})$ as $|\Bx|\rightarrow\infty$.

Note that $p(\Gz,\Gt)= O(|\Gz|^2+|\Gt|^2)$ as $(\Gz,\Gt) \to (0,0)$. Thus, $p(\Bx)=O(|\Bx|^{-2})$ as $|\Bx| \to \infty$, and hence $(\Bh_1, p_1) \in \Mcal$. This completes the proof. \qed

\medskip

It is helpful to write $\Bh_1$ in terms of Cartesian coordinates. By \eqref{eq:bipolar_x_y}, we have
\begin{align*}
    \Psi_1 = A_1 y \Gz + B_1 y \sinh\Gz,
\end{align*}
and hence
\begin{align*}
    \nabla\Psi_1 = A_1 \Gz \Be_y + A_1 y \nabla \Gz + B_1  \sinh\Gz \Be_y + B_1 y  \cosh \Gz \nabla \Gz.
\end{align*}
Then, since $\Bh_1=(\nabla \Psi_1)^\perp$, we have
\beq\label{BhoneCart}
\Bh_1 =   ( A_1 \Gz  + B_1  \sinh\Gz ) \Be_x + (A_1   + B_1   \cosh \Gz) y (\nabla \Gz)^\perp.
\eeq
Here, $(x,y)^\perp = (y,-x)$.

\begin{prop}\label{cor:h1_p1_asymp}
We have
\beq\label{honepone}
\| \Ecal[\Bh_1] \| \lesssim \Gd^{-3/2} \quad\mbox{and}\quad \| p_1\|_\infty \approx \Gd^{-2}
\eeq
as $\Gd \to 0$. In the narrow region $\Pi_\Gd$,
\beq\label{ponenarrow}
\Gs[\Bh_1,p_1](x,y) = -\frac{3}{4} \Gm R \Gd^{-2} \frac{(y^2+3R\Gd)(y^2-R\Gd)}{(y^2+R\Gd)^2} I + O(\Gd^{-3/2}).
\eeq
\end{prop}

\proof
One can see from the explicit forms of the constants $A_1$ and $B_1$ in \eqnref{AB} that
\begin{align}
A_1 &= \frac{3}{8} s^{-3} + O(s^{-1}), \quad B_1 = - \frac{3}{16} s^{-3} + O(s^{-1}).
\label{eq:A1B1_asymp}
\end{align}
Using \eqref{eq:strain_bipolar1}-\eqref{eq:strain_bipolar3} and Proposition \ref{lem:hone} (i), we have
\begin{align}
\Ecal_{\Gz\Gz} &= -h(\Gz,\Gt)(A_1+2B_1 \cosh 2\xi)\cos\Gt,  \label{Ezz} \\
\Ecal_{\Gt\Gt} &= h(\Gz,\Gt)(A_1+2B_1 \cosh 2\xi)\cos\Gt, \label{Ett} \\
  \Ecal_{\Gz\Gt} &= h(\Gz,\Gt) 2B_1 \sinh 2\Gz \sin \Gt. \label{Ezt}
\end{align}

We first estimate $\Ecal_{\Gz\Gz}$. It follows from the Taylor expansions of $\cosh2\Gz$ and $\sinh2\Gz$, and from \eqnref{eq:A1B1_asymp} that
\begin{align*}
    \Ecal_{\Gz\Gz} = - \frac{1+\cos\Gt + O(\Gz^2)}{a} (A_1+ 2B_1 + O(\Gz^2))\cos\Gt .
\end{align*}
Observe from \eqref{eq:A1B1_asymp} that $A_1+2B_1 = O(s^{-1})$. Since $|\Gz| \le s$ and $a,s \approx \sqrt\Gd$,
we have
$$
|\Ecal_{\Gz\Gz}| \lesssim \Gd^{-1}.
$$

Estimates for $\Ecal_{\Gt\Gt}$ and $\Ecal_{\Gz\Gt}$ are simpler. In fact, one can see immediately from \eqnref{Ett} and \eqnref{Ezt} that
$$
|\Ecal_{\Gt\Gt}| = |\Ecal_{\Gz\Gz}| \lesssim \Gd^{-1}
$$
and
$$
|\Ecal_{\Gz\Gt}| \lesssim a^{-1} |B_1 \Gz|  \lesssim \Gd^{-3/2}.
$$
Then \eqnref{Ecalnorm} yields the first estimate in \eqnref{honepone}.

We now consider the pressure $p_1$.
Since $a \approx \sqrt{\Gd}$, we have
$$
|p_1(\Gz,\Gt)| \lesssim \Gd^{-2} (\cosh\Gz \,|\cos\Gt| +1).
$$
Since $|\Gz| \le s \approx \sqrt{\Gd}$ by \eqnref{sGd} if $(\Gz,\Gt) \in D^e$, we have
$$
|p_1(\Gz,\Gt)| \lesssim \Gd^{-2} .
$$
Using the Taylor expansion of $\cosh \Gz$, we see
$$
p_1 = \frac{3}{2} \Gm R \Gd^{-2} \left( \cos \Gt- \frac{1}{2} \cos^2 \Gt \right) +O(\Gd^{-1}).
$$
In particular, we have $\| p_1\|_\infty \gtrsim \Gd^{-2}$, and the second statement in \eqnref{honepone} follows. Now the expansion  \eqnref{ponenarrow} in the narrow region follows from \eqnref{Gsdef} and \eqnref{cosnarrow}.
\qed

\medskip

The expressions for the solution $(\Bh_2,p_2)$ are quite involved even though it is possible to express it explicitly. However, its singular part, which is to be used in the rest of the paper, can be expressed in a rather simple way. To express the singular part, which is denoted by $(\widetilde\Bh_2,\widetilde p_2)$, let
\beq\label{A2C2}
A_2 = -\frac{1}{2s + \sinh 2s}.
\eeq
Then, the components of the velocity field $\widetilde\Bh_2 = \widetilde h_{2\Gz}\Be_\Gz + \widetilde h_{2\Gt} \Be_\Gt$ are given by
\begin{align}
\widetilde h_{2\Gz} &=
A_2 \Gz \beta(\Gz,\Gt),
\label{eq:h2t_Gz}
\\
\widetilde h_{2\Gt} &=A_2\Gz \Ga(\Gz,\Gt) + A_2 \sinh \Gz, \label{eq:h2t_Gt}
\end{align}
and the pressure $\widetilde p_2$  is given by
\beq\label{eq:pressure2}
\widetilde p_2 = -\frac{2\mu}{a} A_2 \sinh \Gz\sin\Gt.
\eeq
Then one can see easily that $(\widetilde\Bh_2,\widetilde p_2)$ belongs to $\Mcal$ and is a solution to the Stokes system. Moreover, $\widetilde\Bh_2$ satisfies
\beq\label{eq:h2_pD_2}
\widetilde{\Bh}_2|_{\p D_i} = \frac{(-1)^i}{2}\Bpsi_2 - C_2 \Bpsi_3,\quad i=1,2,
\eeq
where $\Bpsi_3$ is the one given in \eqnref{Bpsi} and $C_2$ is the constant given by
\beq\label{C2}
C_2 = \frac{\sinh^2 s}{a} A_2.
\eeq
In fact, one can easily check using \eqref{exey} that
\beq\label{eq:Be_Gt_formula}
\Be_\Gt|_{\p D_2} = -\cosh s \Bpsi_2 - \frac{\sinh s}{a}  \Bpsi_3.
\eeq
It then follows from \eqref{eq:h2t_Gz} and \eqref{eq:h2t_Gt} that
\begin{align*}
\widetilde{\Bh}_2|_{\p D_2} &= A_2 s (\Gb \Be_\Gz + \Ga \Be_\Gt) + A_2 \sinh s \Be_\Gt
\\
&=A_2 s (-\Bpsi_2) + A_2 (-\sinh s \cosh s) \Bpsi_2 - \frac{\sinh^2 s}{a} A_2 \Bpsi_3.
\end{align*}
This proves \eqnref{eq:h2_pD_2} on $\p D_2$. \eqnref{eq:h2_pD_2} on $\p D_1$ can be proved in the same way. In Cartesian coordinates, $\widetilde\Bh_2$ is represented in a simple form as
\beq
\widetilde\Bh_2 = -A_2 \Gz \Be_y + A_2 x (\nabla\Gz)^\perp.
\eeq

Some words about how to derive $(\widetilde\Bh_2,\widetilde p_2)$ may be helpful. As in Proposition \ref{lem:hone}, we first derive the relevant stream function $\widetilde\Psi_2$ using the expansion \eqref{eq:general_biharmonic_solution} for the general solution, which turns out to be
\beq\label{eq:stream_singular2}
\widetilde\Psi_2(\Gz,\Gt) = \frac{1}{h(\Gz,\Gt)} A_2 \Gz \sinh \Gz.
\eeq
We then let $(\widetilde\Bh_2,\widetilde p_2)$ be its associated solution to the Stokes system.

Thanks to \eqnref{eq:h2_pD_2}, how to find the solution $(\Bh_2,p_2)$ is clear. Let
$(\Bh_{\mathrm{rot}},p_{\mathrm{rot}})$ be the solution to
\beq\label{eq:def_hrot}
\begin{cases}
\mu \GD \Bh_{\mathrm{rot}} = \nabla p_{\mathrm{rot}} &\quad \mbox{in }D^e,
\\
\nabla \cdot \Bh_{\mathrm{rot}} = 0 &\quad \mbox{in } D^e,
\\[0.3em]
\Bh_{\mathrm{rot}} = \Bpsi_3 &\quad {\p D_1\cup \p D_2},
\\[0.3em]
(\Bh_{\mathrm{rot}},p_{\mathrm{rot}})\in \Mcal.
\end{cases}
\eeq
The existence and uniqueness of the solution are guaranteed by Theorem \ref{thm:ext_diri}. We will construct the stream function for $(\Bh_{\mathrm{rot}},p_{\mathrm{rot}})$ explicitly in subsection \ref{subsec6.1} and prove the following theorem in section \ref{sec:noslip2}.

\begin{theorem}\label{thm:boundedstress_rot}
We have
\beq\label{rotstrain}
\|\Ecal[\Bh_{\mathrm{rot}}]\|_\infty \lesssim 1, \quad \| p_{\mathrm{rot}} \|_\infty \lesssim 1,
\eeq
and
\beq\label{rotstress}
\|\Gs[\mathbf{h}_{\mathrm{rot}},p_{\mathrm{rot}}]\|_\infty \lesssim 1.
\eeq
\end{theorem}

We immediately have the following proposition.

\begin{prop}\label{lem:htwo}
Let $(\widetilde\Bh_2,\widetilde p_2)$ be as given in \eqnref{eq:h2t_Gz}-\eqnref{eq:pressure2} and $C_2$ the constant given in \eqnref{C2}. The solution $(\Bh_2,p_2)$ to \eqnref{Bhj} is given by
\beq
(\Bh_2,p_2)= (\widetilde\Bh_2,\widetilde p_2)+ C_2(\Bh_{\mathrm{rot}},p_{\mathrm{rot}}).
\eeq
\end{prop}

\begin{prop}\label{cor:h2_p2_asymp}
It holds that
\beq\label{htwoptwo}
\| \Ecal[\Bh_2] \|_\infty \approx \Gd^{-1} \quad\mbox{and}\quad \| p_2\|_\infty \approx \Gd^{-1/2},
\eeq
as $\Gd \to 0$. In the narrow region $\Pi_\Gd$,
\beq\label{Gstwonarrow}
\Gs[\Bh_2,p_2](x,y)= \mu \Gd^{-1} \frac{R\Gd}{y^2+R\Gd} \begin{bmatrix}0 & 1 \\ 1 & 0 \end{bmatrix} + O(\Gd^{-1/2}).
\eeq
\end{prop}

\pf
We first note that
\beq\label{Atwo}
A_2= -\frac{1}{4s} + O(s).
\eeq
Since $|\Gz| \le s \approx \sqrt{\Gd}$ and $a \approx \sqrt{\Gd}$, the second estimate in \eqnref{htwoptwo} immediately follows from \eqnref{eq:pressure2}.

Since $a \approx s$ as one can see from \eqnref{s}, it follows from \eqnref{Atwo} and the definition of $C_2$ in \eqnref{C2} that $C_2$ is bounded regardless of $\Gd$. In view of \eqref{rotstrain}, we only need to derive estimates related to $(\widetilde{\Bh}_2, \widetilde p_2)$.
Using \eqref{eq:strain_bipolar1}-\eqref{eq:strain_bipolar3} and \eqnref{eq:stream_singular2},
we have
\begin{align}
\Ecal_{\Gz\Gz}[\widetilde\Bh_2] &= 0, \label{Ezz2}
\\
\Ecal_{\Gt\Gt}[\widetilde\Bh_2] &= 0, \label{Ett2}
\\
\Ecal_{\Gz\Gt}[\widetilde\Bh_2] &=  \frac{\cosh \Gz- \cos\Gt}{a} A_2\cosh \Gz. \label{Ezt2}
\end{align}
We then have from \eqnref{Ezt2} that
$$
|\Ecal_{\Gz\Gt}[\widetilde\Bh_2]| \lesssim \Gd^{-1},
$$
and hence
$$
\| \Ecal[\widetilde\Bh_2] \|_\infty \lesssim \Gd^{-1}.
$$

We see from \eqnref{adef}, \eqnref{sGd}, \eqnref{Atwo} and \eqnref{Ezt2} that
$$
\Ecal_{\Gz\Gt}[\widetilde\Bh_2] = -\frac{1}{4\Gd} (\cosh \Gz- \cos\Gt) + O(1).
$$
In the narrow region $\Pi_\Gd$, we have
$$
\Ecal_{\Gz\Gt}[\widetilde\Bh_2] = -\frac{1}{4\Gd} (1- \cos\Gt) + O(1).
$$
In particular, $|\Ecal_{\Gz\Gt}| \gtrsim \Gd^{-1}$, and the first estimate in \eqnref{htwoptwo} follows.
The asymptotic formula \eqnref{Gstwonarrow} follows from \eqnref{Gsdef}, \eqref{cosnarrow}, \eqnref{Xinarrow} and \eqnref{Gsrel}.
\qed

\section{Proof of Theorem \ref{thm:main1}}\label{sec:1.1}

Thanks to the symmetry of the problem \eqnref{stokes} with $\BU=\BUst=(x,-y)^T$ and $P=0$, the velocity $\Bu$ enjoys the following symmetry:
\begin{align*}
u_{x}(x,y)=u_{x}(x,-y)=-u_{x}(-x,y), \\
u_{y}(x,y)=-u_{y}(x,-y)=u_{y}(-x,y),
\end{align*}
and the pressure $p$ does:
$$
p(x,y) = p(-x,y), \quad p(x,y) = p(x,-y).
$$
Thus, we infer
$$
c_{11}=-c_{21} \quad\mbox{and}\quad c_{i2}=c_{i3}=0 \ \mbox{ for $i=1,2$}.
$$
In other words, we have
\beq\label{eq:u_bdry}
\Bu = - c_{21} \Bpsi_1 \quad\mbox{on } \p D_1, \quad \Bu =  c_{21} \Bpsi_1 \quad\mbox{on } \p D_2.
\eeq
Therefore, the solution $(\Bu,p):=(\Bust, \pst)$ admits the decomposition in terms of the singular function
\beq\label{eq:decomp_simple_u}
\Bu =  \Bv_1 +2  c_{21}{\Bh}_{1}, \quad p =  q_1 +2  c_{21} p_1 \quad \mbox{in } D^e,
\eeq
where $(\Bv_1,q_1)$ is the solution with the no-slip boundary condition, namely,
\beq\label{stokesv}
\begin{cases}
\mu \GD \Bv_1 = \nabla q_1 \quad &\mbox{in }D^e,
\\
\nabla \cdot \Bv_1 = 0 \quad &\mbox{in }D^e,
\\
\Bv_1=0 \quad &\mbox{on } \p D_1 \cup \p D_2, \\
(\Bv_1-\BUst,q_1) \in \Mcal.
\end{cases}
\eeq
 We will construct the stream function for $(\Bv_1,q_1)$ in subsection \ref{subsec6.1} and prove the following theorem in section \ref{sec:noslip2}.

\begin{theorem}\label{thm:boundedstress_v1q1}
Let $(\Bv_1, q_1)$ be the solution to \eqnref{stokesv}. Then, the following estimates hold:
\beq
\|\Ecal[\Bv_1]\|_\infty \lesssim 1, \quad \| q_1 \|_\infty \lesssim 1,
\eeq
and
\beq
\|\Gs[\Bv_1,q_1]\|_\infty \lesssim 1.
\eeq
\end{theorem}

It then follows from \eqref{eq:decomp_simple_u} that
\begin{align}
\Ecal [\Bu] &= 2 c_{21 }\Ecal[\Bh_1] +O(1),\nonumber
\\
p &= 2 c_{21 }p_1 +O(1), \label{eq:pre_asymps}
\\
\Gs[\Bu,p] &= 2 c_{21 }\Gs[\Bh_1,p_1] +O(1), \nonumber
\end{align}
as $\Gd \to 0$. Here, $O(1)$ means that the supremum norms of the remainder terms are bounded on $D^e$ regardless of $\Gd$.
Because of \eqnref{honepone}, it is now sufficient to estimate the constant $c_{21}$.

We first express $c_{21}$ in terms of boundary integrals. To do so, we let
\beq
\Ical_1 := \ds\int_{\p D_2} \Be_x \cdot \Gs[ \Bh_1,p_1]\nu \, dl
\quad\mbox{and}\quad \Jcal_1 :=
\ds \int_{\p D_2} \BU \cdot \Gs [\Bh_1,p_1]\nu \, dl,
\eeq
with $\BU=\BUst=(x,-y)^T$.

\begin{lemma}\label{lem:c21}
We have
\beq
c_{21} = \frac{\Jcal_1}{\Ical_1}.
\eeq
\end{lemma}

\proof
By Green's identity for the Stokes system on $D^e$, we obtain that
\beq\label{1000}
\int_{\p D^e}
(\Bu-\BU)\cdot {\Gs[ \Bh_1,p_1]}  \big|_+ \nu - {\Gs[\Bu-\BU,p]} \big|_+ \nu \cdot\Bh_1
=0.
\eeq
Since $\Bh_1|_{\p D_i} = (-1)^i\frac{1}{2}\Bpsi_1$, it follows from the boundary integral conditions \eqnref{equili} that
$$
\int_{\p D_i} \Gs[\Bu,p] \big|_+ \nu \cdot\Bh_1=0, \quad i=1,2.
$$
Applying Green's identity on $D_i$, we have
$$
\int_{\p D_i} \Gs[\BU,0] \big|_+ \nu \cdot\Bh_1= \int_{\p D_i} \Gs[\BU,p_0] \big|_- \nu \cdot\Bh_1 = 0, \quad i=1,2.
$$
It then follows from \eqnref{1000} that
$$
\int_{\p D^e} (\Bu-\BU)\cdot {\Gs[ \Bh_1,p_1]}  \big|_+ \nu =0,
$$
or equivalently
$$
\int_{\p D^e} \Bu\cdot {\Gs[ \Bh_1,p_1]}  \big|_+ \nu
=\int_{\p D^e} \BU\cdot {\Gs[ \Bh_1,p_1]}  \big|_+ \nu.
$$
By symmetry, we have
$$
\int_{\p D_2} \Bu\cdot {\Gs[ \Bh_1,p_1]}  \big|_+ \nu
=\int_{\p D_2} \BU\cdot {\Gs[ \Bh_1,p_1]}  \big|_+ \nu.
$$
Then the conclusion follows from \eqref{eq:u_bdry}.
\qed

\smallskip

We have the following lemma whose proof is given in Appendix \ref{appendixA}.

\begin{lemma} \label{lem:asymp_I1_J1}
As $\Gd\rightarrow 0$, we have
\beq\label{Icalone}
\Ical_1 = -\frac{3\pi \mu}{2} \left(\frac{R}{\Gd}\right)^{3/2} + O(\Gd^{-1/2}),
\eeq
and
\beq\label{Jcalone}
\Jcal_1 = -3\pi \mu R + O(\Gd).
\eeq
\end{lemma}

As an immediate consequence of Lemmas \ref{lem:c21} and \ref{lem:asymp_I1_J1}, we have the following corollary:

\begin{cor}\label{cor:c21_asymp}
As $\Gd\rightarrow 0$, we have
\beq\label{c21}
c_{21} = \frac{2}{\sqrt{R}} \Gd^{3/2} + O(\Gd^{5/2}).
\eeq
\end{cor}

Now Theorem \ref{thm:main1} follows from Proposition \ref{cor:h1_p1_asymp}, \eqnref{eq:pre_asymps}, and Corollary \ref{cor:c21_asymp}.
\medskip

\section{Proof of Theorem \ref{thm:main2}}\label{sec:1.2}

Assume that $\BU(x,y)=\BUsh=(y,x)^T$. We write $(\Bu,p)$ for $(\Bush,\psh)$ for ease of notation. In this case the velocity $\Bu$ satisfies
\begin{align}
u_{x}(x,y)=-u_{x}(x,-y)=u_{x}(-x,y),\nonumber \\
u_{y}(x,y)=u_{y}(x,-y)=-u_{y}(-x,y),
\end{align}
and the pressure $p$ satisfies:
$$
p(x,y) = -p(-x,y), \quad p(x,y) = -p(x,-y).
$$
Then, we see that $c_{22}=-c_{12}$, $c_{23}=c_{13}$ and $c_{i1}=c_{i1}=0$ for $i=1,2$. As a result, we have ,
\beq\label{eq:u_bdry2}
\Bu = - c_{22} \Bpsi_2 + c_{23} \Bpsi_3 \quad\mbox{on } \p D_1, \quad \Bu = c_{22} \Bpsi_2 + c_{23}\Bpsi_3  \quad\mbox{on } \p D_2.
\eeq

Let us decompose the solution $(\Bu, p)$ in $D^e$ as
\beq\label{eq:decomp_simple_u_two}
(\Bu, p)  =  (\Bv_2,q_2) + 2 c_{22}({\Bh}_{2},p_2) + c_{23}(\Bh_{\mathrm{rot}},p_{\mathrm{rot}}),
\eeq
where $(\Bv_2,q_2)$ is the solution to
\beq\label{Bv2}
\begin{cases}
\mu \GD \Bv_2 = \nabla {q}_2 &\quad \mbox{in }D^e,
\\
\nabla \cdot \Bv_2 = 0 &\quad \mbox{in } D^e,
\\
\Bv_2 = 0 &\quad \p D_1 \cup \p D_2,
\\
(\Bv_2-\BUsh,q_2) \in \Mcal,
\end{cases}
\eeq
and $(\Bh_{\mathrm{rot}},q_{\mathrm{rot}})$ is the solution to \eqnref{eq:def_hrot}. Note that $\Bv_2$ also satisfies the no-slip boundary condition like $\Bv_1$. We will construct the stream function for $(\Bv_2,q_2)$ together with those for $(\Bh_{\mathrm{rot}},q_{\mathrm{rot}})$ and $(\Bv_1,q_1)$ in subsection \ref{subsec6.1} and prove the following theorem in section \ref{sec:noslip2}.

\begin{theorem}\label{thm:boundedstress_v2q2}
We have
\beq
\|\Ecal[\Bv_2]\|_\infty \lesssim 1, \quad \| q_2 \|_\infty \lesssim 1,
\eeq
and
\beq
\|\Gs[\Bv_2,q_2]\|_\infty \lesssim 1.
\eeq
\end{theorem}

It follows from \eqref{eq:decomp_simple_u_two} that
\begin{align}
\Ecal [\Bu] &= 2 c_{22 }\Ecal[\Bh_2] + c_{23 }\Ecal[\Bh_{\mathrm{rot}}] +O(1),\nonumber
\\
p &= 2 c_{22 }p_2 + c_{23}p_{\mathrm{rot}} +O(1), \label{eq:pre_asymps2}
\\
\Gs[\Bu,p] &= 2 c_{22 }\Gs[\Bh_2,p_2]+c_{23 }\Gs[\Bh_{\mathrm{rot}},p_{\mathrm{rot}}] +O(1), \nonumber
\end{align}
as $\Gd\rightarrow 0$.

As before, we represent the constant $c_{22}$ using the integrals
\begin{align}
\Ical_{2j} &:= \ds\int_{\p D_2} \Bpsi_j \cdot \Gs[ \Bh_2,p_2]\nu \, dl,\quad j=2,3,
\\
\Ical_{\mathrm{rot}} &:= \ds\int_{\p D_2} \Bpsi_3 \cdot \Gs[ \Bh_{\mathrm{rot}},p_{\mathrm{rot}}]\nu \, dl,
\\  \Jcal_2 &:=
\ds \int_{\p D_2} \BU \cdot \Gs [\Bh_2,p_2]\nu \, dl,
\\  \Jcal_{\mathrm{rot}} &:=
\ds \int_{\p D_2} \BU \cdot \Gs [\Bh_{\mathrm{rot}},p_{\mathrm{rot}}]\nu \, dl, \label{ItwoJtwo}
\end{align}
where $\BU=\BUsh=(y,x)^T$.
We have the following lemma whose proof is similar to the one of Lemma \ref{lem:c21}.
\begin{lemma}
We have
$$
 \begin{bmatrix}
 \Ical_{22} & \Ical_{23}
 \\
 \Ical_{23} & \Ical_{\mathrm{rot}}
 \end{bmatrix}
 \begin{bmatrix}
 c_{22}
 \\
 c_{23}
 \end{bmatrix}
 =
 \begin{bmatrix}
 \Jcal_2\\ \Jcal_{\mathrm{rot}}
 \end{bmatrix}.
$$
\end{lemma}
\proof
As in the proof of Lemma \ref{lem:c21}, we have
$$
\int_{\p D_2} \Bu\cdot {\Gs[ \Bh_2,p_2]}  \big|_+ \nu
=\int_{\p D_2} \BU\cdot {\Gs[ \Bh_2,p_2]}  \big|_+ \nu,
$$
and
$$
\int_{\p D_2} \Bu\cdot {\Gs[ \Bh_{\mathrm{rot}},p_{\mathrm{rot}}]}  \big|_+ \nu
=\int_{\p D_2} \BU\cdot {\Gs[ \Bh_{\mathrm{rot}},p_{\mathrm{rot}}]}  \big|_+ \nu.
$$
Then, by \eqref{eq:u_bdry2}, we see
$$
\Ical_{22} c_{22} + \Ical_{23} c_{23} = \Jcal_2,
$$
and
$$
\int_{\p D_2} \Bpsi_2\cdot {\Gs[ \Bh_{\mathrm{rot}},p_{\mathrm{rot}}]}  \big|_+ \nu \cdot c_{22} + \Ical_{\mathrm{rot}} c_{23} = \Jcal_{\mathrm{rot}}.
$$
Then, Green's identity yields
\begin{align}
\int_{\p D_2} \Bpsi_2\cdot {\Gs[ \Bh_{\mathrm{rot}},p_{\mathrm{rot}}]}  \big|_+ \nu &=
\int_{\p D^e} \Bh_2\cdot {\Gs[ \Bh_{\mathrm{rot}},p_{\mathrm{rot}}]}  \big|_+ \nu =
\int_{\p D^e} \Bh_{\mathrm{rot}}\cdot {\Gs[ \Bh_{2},p_{2}]}  \big|_+ \nu \nonumber \\
&=
\int_{\p D_2} \Bpsi_3\cdot {\Gs[ \Bh_{2},p_{2}]}  \big|_+ \nu = \Ical_{23},
\label{eq:I23_another}
\end{align}
and hence the conclusion follows.
\qed

We have the following lemma whose proof is given in Appendix \ref{sec:appendixB}.
Let
\begin{align}
&f_0(x) : =  \frac{ 4e^{-x} \sinh^2 x (\cosh x + \sinh x) -4x^2}{x^3(\sinh 2x + 2x)}, \label{def_small_f}
\\
&g_0(x) :=\frac{4x}{\sinh 2x + 2x}, \label{def_small_g}
\end{align}
and let
\beq\label{F0G0}
F_0:= \int_0^\infty f_0(x) dx, \quad G_0:=\int_0^\infty g_0(x) dx.
\eeq

\begin{lemma}\label{lem:asymp_I2_J2}
As $\delta\rightarrow 0$, we have
\begin{align}
&\Ical_{22} =   -\pi\mu\sqrt{\frac{R}{\Gd}}+O(1),
\label{eq:I22}
\\
&\Ical_{23} =  \frac{\pi\mu R}{F_0} + O(\sqrt{\Gd}),
\label{eq:I23}
\\
&\Ical_{\mathrm{rot}} = -\frac{4\pi\mu R^2}{F_0} + O(\sqrt{\Gd}).
\label{eq:Irot}
\\
&\Jcal_2  = -\pi \mu R \left( 1-\frac{1-G_0}{F_0}\right) +O(\sqrt{\Gd}), \label{eq:J2}
\\
&\Jcal_{\mathrm{rot}} = - 4\pi\mu R^2 \frac{1-G_0}{F_0} + O(\sqrt{\Gd}).
\label{eq:Jrot}
\end{align}
\end{lemma}

\smallskip

As an immediate consequence, the following corollary holds:

\begin{cor}\label{cor:c22_asymp}
As $\Gd\rightarrow 0$, we have
\beq\label{c22}
c_{22} =    \sqrt{R\Gd} +O(\Gd),
\quad
c_{23} = O(1).
\eeq
\end{cor}

Now, Theorem \ref{thm:main2} follows from Theorem \ref{thm:boundedstress_rot}, Proposition \ref{cor:h2_p2_asymp},  \eqref{eq:pre_asymps2} and Corollary \ref{cor:c22_asymp}.
One can also see that the decomposition formula \eqnref{Bpcomp} for the solution $(\Bu,p)$ is an immediate consequence of \eqnref{eq:pre_asymps}, \eqnref{c21}, \eqnref{eq:pre_asymps2} and \eqnref{c22}.

\section{No blow-up with no-slip boundary conditions I}\label{sec:noslip1}

In this and next sections, we show that the stress tensor does not blow up under the no-slip boundary condition, that is, we prove Theorems \ref{thm:boundedstress_rot}, \ref{thm:boundedstress_v1q1} and \ref{thm:boundedstress_v2q2}.  Theorem \ref{thm:boundedstress_rot} is for the problem with the boundary condition given by $\Bpsi_3$, and Theorems \ref{thm:boundedstress_v1q1} and \ref{thm:boundedstress_v2q2} for those with the no-slip boundary conditions. For doing so we first construct solutions $(\Bv_j,q_j)$, $j=1,2$, and $(\Bh_{\mathrm{rot}},p_{\mathrm{rot}})$ by using the stream function formulation and bipolar coordinates. To avoid notational confusion, we denote the stream functions by $\Phi$ in this section instead of $\Psi$ which was used in previous sections.

\subsection{Construction of stream functions}\label{subsec6.1}

In the following three lemmas we present stream functions for $(\Bv_1-\BUst,q_1)$, $(\Bv_2-\BUsh,q_2)$, and $(\Bh_{\mathrm{rot}},p_{\mathrm{rot}})$. Each stream function is found using the general form  \eqref{eq:general_biharmonic_solution} and matching the boundary conditions using the formula \eqref{eq_velo1} and \eqref{eq_velo2} for $\Bu$ components of the solution.

\subsubsection{Stream function for $(\Bv_1-\BUst,q_1)$}

\begin{lemma}\label{lem:Phi_1}
Let $\Phi_1$ be the stream function associated with the solution $(\Bv_1-\BUst,q_1)$.
We have
\begin{align}
 (h\Phi_1)(\Gz,\Gt) &=
 a_1  \sinh 2\Gz\sin\Gt + b_1 \Gz\sin\Gt  \nonumber
\\
&\quad +\sum_{n=2}^\infty \Big({a}_n \sinh (n+1)\Gz+ {b}_n \sinh(n-1)\Gz \Big)\sin n\Gt,
\label{eq:stream_v1}
\end{align}
where
\begin{align*}
    a_1 &=-\frac{2ae^{-s}(\sinh s -s e^{-s})}{\sinh 2s - 2s \cosh 2s},
    \\
    b_1 &=\frac{4a \sinh^2 s}{\sinh 2s - 2s \cosh 2s},
    \\
    a_n &=- \frac{2a(e^{-ns} \sinh ns - e^{-s} n \sinh s)}{\sinh 2ns-n\sinh 2s}, \quad n\geq 2,
    \\
    b_n &= \frac{2a(e^{-ns} \sinh ns - e^{s} n \sinh s)}{\sinh 2ns-n\sinh 2s}, \quad n\geq 2.
\end{align*}
\end{lemma}

\proof
We need to show that the solution $(\Bv_1, q_1)$ constructed from $\Phi_1$ satisfies the no-slip condition $\Bv_1=0$ on $\p D_1$ and $\p D_2$, and $(\Bv_1-\BUst,q_1) \in \Mcal$.

We first observe that $\Phi^{0}_{1} := xy$ is the stream function associated to the background solution $(\BU,0)$. Here and throughout this proof $\BU=\BUst$. In fact, $\BU=(x,-y) = (\nabla \Phi_1^0)^\perp$ and a harmonic conjugate of $\mu\GD \Phi^{0}_{1}=0$ is constant. We see from \eqref{eq:bipolar_x_y} that
$$
\Phi^{0}_{1} = \frac{a^2 \sinh \Gz \sin \Gt}{(\cosh\Gz-\cos\Gt)^2}
$$
in bipolar coordinates.

Notice that $\Phi_1^0$ has the odd symmetry in both $\Gz$ and $\Gt$. So we look for $\Phi_1$ with the same symmetric property. We assume that $\Phi_1$ of the form \eqnref{eq:stream_v1} which is the part with such symmetry of the general solution \eqref{eq:general_biharmonic_solution}, and determine the coefficients $a_n$ and $b_n$ from the no-slip boundary condition.

For that, define $\Phi_1^{\mathrm{tot}}$ by
$$
\Phi_1^{\mathrm{tot}} := \Phi_1^0 + \Phi_1,
$$
so that  $\Phi_1^{\mathrm{tot}}$ is the stream function  associated with $(\Bv_1,q_1)$. If we write $\Bv_1 = v_{1\Gz} \Be_\Gz + v_{1\Gt} \Be_\Gt$, then the no-slip boundary condition becomes
\begin{align}
v_{1\Gz} = 0 \quad \mbox{on }\Gz=\pm s,\label{eq:v1Gz_zero_bd}
\\
 \quad v_{1\Gt}=0 \quad \mbox{on }\Gz=\pm s.\label{eq:v1Gt_zero_bd}
\end{align}
Then, from the formula \eqref{eq:tan_deri_bipolar} for the tangential derivative and \eqref{eq_velo1} for the stream function in bipolar coordinates, we have
\begin{align*}
0={v}_{1\Gz}|_{\Gz =\pm s} &=  -( h \p_\Gt \Phi_{1}^{\mathrm{tot}})|_{\Gz=\pm s}
 =\mp \p_T \Phi_{1}^{\mathrm{tot}}\Big|_{\Gz=\pm  s}.
\end{align*}
This amounts to $\Phi_{1}^{\mathrm{tot}}$ being constant on $\{\Gz=s\}$ and $\{\Gz=-s\}$. Since $\Phi_{1}^{\mathrm{tot}}$ is  odd in $\Gz$, we further require that
$$
\Phi_{1}^{\mathrm{tot}} = 0 \quad \mbox{on } \Gz=\pm s.
$$
We also have from \eqref{eq_velo2} and \eqref{eq:v1Gt_zero_bd} that on $\{\Gz=\pm s\}$
$$
0=	{v}_{1\Gt} = h \p_\Gz \Phi_{1}^{\mathrm{tot}} = \left( \p_\Gz - \frac{\sinh\Gz}{\cosh\Gz-\cos\Gt}\right)(h\Phi_{1}^{\mathrm{tot}})
=\p_{\Gz} (h\Phi_{1}^{\mathrm{tot}}).
$$
Thus the no-slip boundary condition is fulfilled if
$$
\begin{cases}
	h\Phi_{1}^{\mathrm{tot}} = 0 &\quad \mbox{on } \Gz=\pm s,
	\\
	\p_{\Gz} (h\Phi_{1}^{\mathrm{tot}}) = 0&\quad \mbox{on } \Gz=\pm s,
\end{cases}
$$
or equivalently
\beq\label{noslipv1}
\begin{cases}
	h\Phi_{1} = - h\Phi_1^0 &\quad \mbox{on } \Gz=\pm s,
	\\
	\p_{\Gz} (h\Phi_{1}) = - \p_\Gz(h\Phi_1^0)&\quad \mbox{on } \Gz=\pm s.
	\end{cases}
\eeq

Note that
\begin{align}
(h \Phi_1^0)(\Gz,\Gt) &= \frac{a \sinh\Gz\sin\Gt}{\cosh \Gz-\cos\Gt}
=2a \sinh\Gz \sum_{n=1}^\infty e^{-n |\Gz|} \sin n\Gt, \quad \Gz\neq 0.
\end{align}
We then see from \eqnref{eq:stream_v1} that \eqnref{noslipv1} is equivalent to the following linear systems for $a_n$ and $b_n$:
\begin{align}
\begin{bmatrix}
\sinh 2s && s
\\
2\cosh 2s && 1	
\end{bmatrix}
\begin{bmatrix}
a_1
\\
b_1	
\end{bmatrix}
=
\begin{bmatrix}
-2a \sinh s e^{-s}
\\
-2a \cosh s e^{-s} + 2a \sinh s e^{-s}	
\end{bmatrix},
\end{align}
and
\begin{align}
&\begin{bmatrix}
\sinh (n+1)s & \sinh(n-1)s
\\
(n+1)\cosh (n+1)s & (n-1)\cosh (n-1)s	
\end{bmatrix}
\begin{bmatrix}
a_n
\\
b_n	
\end{bmatrix}
\nonumber
\\
&=
\begin{bmatrix}
-2a \sinh s e^{-ns}
\\
-2a \cosh s e^{-ns} + 2a \sinh s n e^{-ns}	
\end{bmatrix}, \quad n \ge 2.
\end{align}
Solving these linear systems yield the expressions for $a_n$ and $b_n$.

We now show
\beq\label{voneMcal}
(\Bv_1-\BU,q_1) \in \Mcal.
\eeq
We first prove
\beq\label{eq:claim_uo_bdd}
(\Bv_1-\BU)(\Bx)=O(|\Bx|^{-1}), \quad |\Bx|\rightarrow\infty.
\eeq
Since $|\Bx|\rightarrow\infty$ is equivalent to $(\Gz,\Gt)\rightarrow(0,0)$, it is equivalent to proving
\beq
(\Bv_1-\BU)(\Gz,\Gt) = O(|\Gz|+|\Gt|), \quad (\Gz,\Gt)\rightarrow(0,0).
\eeq

One can see from the explicit forms of $a_n$ and $b_n$ that there is a constant $C$ independent of $n$ ($C$ may  dependent on $s$) such that
$$
|a_n| + |b_n| \le C n e^{-2ns}
$$
for all $n$. Thus for any positive number $k$ there is a constant $C$ such that
\beq\label{eq:anbnsum_conv}
\sum_{n=1}^\infty n^k e^{ns}(|a_n| + |b_n|) \le C.
\eeq
The constant $C$ may differ at each appearance.

If we write $\Bv_1-\BU = f_{\Gz}\Be_\Gz + f_{\Gt}\Be_\Gt$, then it follows from \eqref{eq_velo1} and \eqref{eq_velo2} that
\begin{align}
f_{\Gz} &=  - h \p_\Gt \Phi_1 =\left(- \p_\Gt + F \right)(h\Phi_1), \label{eq:uo_Gz} \\
f_{\Gt} &= +h \p_\Gz \Phi_1 = \left( \p_\Gz - G \right)(h\Phi_1), \label{eq:uo_Gt}
\end{align}
where
\beq
F:= \frac{\sin\Gt}{\cosh\Gz-\cos\Gt}, \quad G:= \frac{\sinh\Gz}{\cosh\Gz-\cos\Gt}.
\eeq
According to \eqnref{eq:stream_v1}, $h\Phi_1$ can be written as
\beq
(h \Phi_1)(\Gz,\Gt) = a_1  \sinh 2\Gz\sin\Gt + b_1 \Gz\sin\Gt
+\sum_{n=2}^\infty \left( {a}_n w_n^+(\Gz,\Gt) + {b}_n w_n^-(\Gz,\Gt) \right) ,	
\eeq
where
\beq
w_n^\pm(\Gz,\Gt):= \sinh (n\pm 1)\Gz \sin n\Gt.
\eeq

One can see that
\beq\label{eq:def_wnpm}
|w_n^\pm(\Gz,\Gt)| \lesssim n^2e^{n s} |\Gz\Gt|	.
\eeq
It thus follows from \eqnref{eq:anbnsum_conv} that
\beq\label{eq:hPhio_estim}
|(h\Phi_1)(\Gz,\Gt)| \lesssim |\Gz\Gt| + |\Gz\Gt| \sum_{n=2}^\infty n^2 e^{ns}(|a_n|+|b_n|) \lesssim |\Gz\Gt|.
\eeq
Similarly, one can show that there is $C$ independent of $(\Gz,\Gt)$ such that
\begin{align}
&|\p_\Gz (h\Phi_1)| \le C |\Gt|,
\quad
|\p_\Gt (h\Phi_1)| \le C  |\Gz|, \nonumber
\\
&|\p_\Gz^2 (h\Phi_1)| \le C |\Gz\Gt|,
\quad
|\p_\Gt^2 (h\Phi_1)| \le C  |\Gz\Gt|,
\quad |\p_\Gz\p_\Gt (h\Phi_1)| \le C .
\label{eq:hPhio_d_estim}
\end{align}

Since
\beq\label{FGest}
\left| F \right| \approx  \frac{|\Gt|}{\Gz^2+\Gt^2}, \quad \left| G \right| \approx \frac{|\Gz|}{\Gz^2+\Gt^2}
\eeq
as $(\Gz,\Gt) \to 0$, we have from \eqref{eq:uo_Gz}, \eqref{eq:uo_Gt}, \eqnref{eq:hPhio_estim} and \eqnref{eq:hPhio_d_estim} that
\beq\label{eq:uoGz_estim}
|f_{\Gz}| + |f_\Gt| \le C (|\Gz| + |\Gt|)
\eeq
for some constant $C$ (depending on $s$, and hence on $\Gd$), which implies \eqref{eq:claim_uo_bdd}.

Next, we prove
\beq\label{eq:claim_grad_uo_decay}
\nabla (\Bv_1-\BU)(\Bx) = O(|\Bx|^{-2}), \quad |\Bx|\rightarrow\infty,
\eeq
or equivalently
\beq
\nabla(\Bv_1-\BU)(\Gz,\Gt) = O(\Gz^2+\Gt^2), \quad (\Gz,\Gt)\rightarrow(0,0).
\eeq
Since $\Bv_1-\BU = f_{\Gz}\Be_\Gz + f_{\Gt}\Be_\Gt$, we have
$$
|\nabla(\Bv_1-\BU)| \leq C(| \nabla f_\Gz| + |f_\Gz \nabla \Be_\Gz| + |\nabla f_\Gt|+|f_\Gt\nabla\Be_\Gt|).
$$

Lemma \ref{lem:grad_eGz_eGt_estim} and \eqref{eq:uoGz_estim} yield
$$
|\nabla\Bv_1-\BU| \lesssim | \nabla f_\Gz| + |\nabla f_\Gt| + (\Gz^2+\Gt^2).
$$
We see from \eqref{eq:uo_Gz}
$$
\p_\Gz f_{\Gz} =\left(- \p_\Gz\p_\Gt + F \p_\Gz + \p_\Gz F \right)(h\Phi_1).
$$
One can see easily that $|\p_\Gz F| \lesssim (\Gz^2+\Gt^2)^{-1}$. Thus we obtain from \eqref{eq:hPhio_estim}, \eqref{eq:hPhio_d_estim} and \eqnref{FGest}
$$
	\p_\Gz f_{\Gz}= O(1) .
$$
Similarly, one can show
\begin{align*}
	\p_\Gt f_{\Gz}= O(1), \quad \p_\Gz f_{\Gt}= O(1), \quad \p_\Gt f_{\Gt}= O(1).
\end{align*}
Therefore, we have
\begin{align}
\nabla f_\Gz= O(|h \p_\Gz f_\Gz| + |h \p_\Gt f_\Gz|) = O(|h|) = O(\Gz^2+\Gt^2),
\\
\nabla f_\Gt= O(|h \p_\Gz f_\Gt| + |h \p_\Gt f_\Gt|) = O(|h|) = O(\Gz^2+\Gt^2).
\end{align}
This proves \eqref{eq:claim_grad_uo_decay}.

We now prove the estimate of the pressure:
\beq\label{eq:claim_p_uo_decay}
q_1(\Bx) = O(|\Bx|^{-2}), \quad |\Bx|\rightarrow\infty,	
\eeq
or equivalently,
\beq
q_1 = O(\Gz^2+\Gt^2), \quad (\Gz,\Gt)\rightarrow(0,0).
\eeq

Let
\beq\label{wn}
w_n(\Gz,\Gt):=\sinh n\Gz \sin n \Gt, \quad
\tilde{w}_n(\Gz,\Gt)=\cosh n\Gz \cos n \Gt.
\eeq
The pressure $q_1$ is given by
\begin{align}
q_1 &=C -a_1\frac{2\mu}{a} (2 \tilde{w}_1-\tilde{w}_2 ) +b_1\frac{2\mu}{a}\tilde{w}_1
 - \frac{2\mu}{a}\sum_{n=2}^\infty((n+1)a_n - (n-1)b_n )\tilde{w}_n \nonumber
    \\
 &\quad + \frac{2\mu}{a}\sum_{n=2}^\infty n a_n \tilde{w}_{n+1}
     - \frac{2\mu}{a}\sum_{n=2}^\infty n b_n \tilde{w}_{n-1}, \label{qone}
\end{align}
for some constant $C$. In fact, one can see from \eqref{eq:Laplacian_Psi} that
\begin{align*}
    \GD \Phi_1 &= a_1\frac{2}{a} (2w_1-w_2) -b_1\frac{2}{a}w_1
    + \frac{2}{a}\sum_{n=2}^\infty((n+1)a_n - (n-1)b_n )w_n
    \\&\quad - \frac{2}{a}\sum_{n=2}^\infty n a_n w_{n+1}
     + \frac{2}{a}\sum_{n=2}^\infty n b_n w_{n-1}.
\end{align*}
Since $q_1$ is a harmonic conjugate of $\mu\GD\Phi_1$ and $-\tilde{w}_n$ is a harmonic conjugate of $w_n$, \eqnref{qone} follows.

We choose the constant $C$ to be
\beq\label{qonepinfty}
C = \frac{2\mu}{a}a_1 - \frac{2\mu}{a}b_1+\frac{2\mu}{a}\sum_{n=2}^\infty (a_n+b_n).
\eeq
Then, $q_1$ take the form
\begin{align}
q_1 &=C-a_1\frac{2\mu}{a} (2 (\tilde{w}_1-1)-(\tilde{w}_2-1) ) +b_1\frac{2\mu}{a}(\tilde{w}_1-1)\nonumber
 \\
 &\quad - \frac{2\mu}{a}\sum_{n=2}^\infty((n+1)a_n - (n-1)b_n )(\tilde{w}_n -1)
 \nonumber   \\&\quad
    +     \frac{2\mu}{a}\sum_{n=2}^\infty n a_n (\tilde{w}_{n+1}-1)
     - \frac{2\mu}{a}\sum_{n=2}^\infty n b_n (\tilde{w}_{n-1}-1) . \label{qone2}
\end{align}

Note that
\beq\label{eq:wtn_estim}
|\tilde{w}_n(\Gz,\Gt) - 1| \lesssim n^2 e^{n s}(\Gz^2 + \Gt^2).
\eeq
This together with \eqnref{eq:anbnsum_conv} yields
\begin{align*}
|q_1|&\lesssim \Big( 1 + \sum_{n=2}^\infty  n^3 e^{ns}(|a_n|+|b_n|)\Big)|\Gz^2+\Gt^2 | =O(\Gz^2+\Gt^2).
\end{align*}
This proves \eqref{eq:claim_p_uo_decay} and hence \eqnref{voneMcal}. The proof is completed.
\qed

\subsubsection{Stream function for $(\Bv_2-\BUsh,q_2)$}

\begin{lemma}\label{lem:Phi_2}
The stream function $\Phi_2$ associated with the solution $(\Bv_2-\BUsh,q_2)$ is given by
\begin{align}
(h\Phi_2)(\Gz,\Gt)
&= K_v(\cosh\Gz -\cos\Gt)\ln (2\cosh \Gz - 2\cos \Gt) +c_0 \cosh \Gz + d_0 \Gz \sinh \Gz \nonumber
\\
&\quad
+\sum_{n=1}^\infty \big(c_n \cosh (n+1) \Gz + d_n \cosh (n-1)\Gz \big)\cos n\Gt ,
\label{eq:W_biharmonic_solution}
\end{align}
where
\begin{align*}
K_v=\frac{a(1-\tanh s -\frac{2\sinh^2 s}{2s+ \sinh 2s}-M')}{\frac{1}{2} + \frac{s(\sinh 2s -2 \tanh s)}{2s +\sinh 2s} + M}
\end{align*}
with
\begin{equation} \label{M}
    M=
     \sum_{n=2}^\infty \frac{4n \sinh s \cosh s + e^{-n s}\sinh n s - 4n^2 \sinh^2 s }{n(n^2-1)(\sinh 2n s + n \sinh 2s)},
     \end{equation}
     and \begin{equation} \label{M'}
    M'=\sum_{n=2}^\infty \frac{4 n \sinh^2 s }{\sinh 2ns + n \sinh 2s},
    \end{equation}
\begin{align*}
    c_0 &=-\frac{a}{2}+\frac{a \sinh^2 s}{\sinh s \cosh s+s} + K_v\frac{-1+e^{-2s} -2s(1+s)}{2s + \sinh 2s}, \\
    d_0 &= \frac{a}{\sinh s \cosh s +s} -K_v \frac{\sinh^2 s}{s+\cosh s \sinh s},
    \\
    c_1 &=a(-1+\coth 2s) + K_v \frac{1}{1+e^{2s}},
    \\
    d_1 &= \frac{a}{2}-\frac{a}{\sinh 2s} + K_v(1+s-\frac{\tanh s}{2}),
    \\
    c_n &=\frac{2a(e^{-ns}\cosh ns- e^{-s}n \sinh s)}{\sinh 2ns + n \sinh 2s}+2K_v \frac{e^{-ns} \sinh n s + e^{-s}n\sinh s }{n(n+1)(\sinh 2ns + n \sinh 2s)},
    \\
    d_n &=-\frac{2a(e^{-ns}\cosh ns- e^{s}n \sinh s)}{\sinh 2ns + n \sinh 2s} - 2K_v\frac{e^{-n s} \sinh ns + e^s n \sinh s}{n(n-1) (\sinh 2ns + n \sinh 2s)}.
\end{align*}
\end{lemma}

\proof
Like the proof of Lemma \ref{lem:Phi_1}, one can see that the stream function associated with the background solution $(\BU,0)$ is given by
\beq
 \Phi^{0}_{2} = \frac{1}{2}(-x^2+y^2).
\eeq
Here and throughout this proof, $\BU=\BUsh$. One can see from \eqref{eq:bipolar_x_y} that
\beq\label{eq:Phi20_bipolar}
\Phi^{0}_{2} = \frac{1}{2}\frac{a^2 (-\sinh^2 \Gz+ \sin^2 \Gt)}{(\cosh\Gz-\cos\Gt)^2}.
\eeq

Since $\Phi_2^0$ has the even symmetry in both $\Gz$ and $\Gt$, we seek $\Phi_2$ in the form \eqnref{eq:W_biharmonic_solution} which has the same symmetric property. Let
\begin{align}
 (h\Phi_{2}^K)(\Gz,\Gt) &:= K_v (\cosh\Gz -\cos\Gt)\ln (2\cosh \Gz - 2\cos \Gt),\label{eq:Phi2K}\\
 (h\Phi_{2}^F)(\Gz,\Gt) &:=  c_0 \cosh \Gz + d_0 \Gz \sinh \Gz
\nonumber
\\
&\quad
+\sum_{n=1}^\infty \big(c_n \cosh (n+1) \Gz + d_n \cosh (n-1)\Gz \big)\cos n\Gt,
\label{eq:Phi2F_fourier}
\end{align}
so that
\begin{align}
& \Phi_{2} = \Phi_2^K + \Phi_2^F. \label{eq:Phi2_fourier}
\end{align}

Let
$$
\Phi_2^{\mathrm{tot}} := \Phi_2^0 + \Phi_2.
$$
Then, $\Phi_2^{\mathrm{tot}}$ is the stream function  associated with $(\Bv_2,q_2)$.
We determine the coefficients $c_n$ and $d_n$ from the no-slip boundary condition $\Bv_2=0$ on $\p D_1$ and $\p D_2$. One can show as in the proof of Lemma \ref{lem:Phi_1} that this condition is fulfilled if
\beq
\begin{cases}
	h\Phi_{2}^{\mathrm{tot}} = 0 &\quad \mbox{on } \Gz=\pm s,
	\\
	\p_{\Gz} (h\Phi_{2}^{\mathrm{tot}}) = 0&\quad \mbox{on } \Gz=\pm s.
	\end{cases}
\eeq
In other words,
\beq\label{systemv2}
\begin{cases}
	h\Phi_{2}^F = - h\Phi_2^0 - h\Phi_2^K &\quad \mbox{on } \Gz=\pm s,
	\\
	\p_{\Gz} (h\Phi_{2}^F) = - \p_\Gz(h\Phi_2^0) - \p_\Gz(h\Phi_2^K)&\quad \mbox{on } \Gz=\pm s.
	\end{cases}
\eeq

Let
\begin{align}
(h \Phi_2^0)(\Gz,\Gt) &= \frac{1}{2}\frac{a^2 (-\sinh^2\Gz+\sin^2\Gt)}{\cosh \Gz-\cos\Gt}\nonumber
\\
&=\frac{a}{2} e^{-|\Gz|} + \frac{a}{2}e^{-2|\Gz|}\cos\Gt - a\sinh|\Gz|\sum_{n=2}^\infty e^{-n|\Gz|}\cos n\Gt\nonumber
\\
&=:\sum_{n=0}^\infty \phi_{n}^0(\Gz) \cos n\Gt,
\end{align}
and
\begin{align}
(h \Phi_2^K)(\Gz,\Gt)
&=K_v(|\Gz|\cosh\Gz + e^{-|\Gz|}) -K_v\Big(1+\frac{e^{-2|\Gz|}}{2}+|\Gz|\Big)\cos\Gt
\nonumber\\
&\quad + K_v\sum_{n=2}^\infty \Big(\frac{e^{-(n-1)|\Gz|}}{n-1}-2\cosh \Gz \frac{e^{-n|\Gz|}}{n} +\frac{e^{-(n+1)|\Gz|}}{n+1}\Big)\cos n\Gt\nonumber
\\
&=:\sum_{n=0}^\infty \phi_{n}^K(\Gz) \cos n\Gt.
\end{align}
Then one can infer from \eqnref{systemv2} that the following system of equations for $c_n, d_n$ hold:
\begin{align*}
\begin{bmatrix}
\cosh s && s\sinh s
\\
\sinh s && \sinh s + s\cosh s	
\end{bmatrix}
\begin{bmatrix}
c_0
\\
d_0	
\end{bmatrix}
&=
\begin{bmatrix}
-\phi_{0}^K(s) -\phi_{0}^K(s)
\\
-(\phi_{0}^K)'(s) -(\phi_{0}^K)'(s)
\end{bmatrix},
\\
\begin{bmatrix}
\cosh 2s && 1
\\
2\sinh 2s && 0
\end{bmatrix}
\begin{bmatrix}
c_1
\\
d_1	
\end{bmatrix}
&=
\begin{bmatrix}
-\phi_{1}^K(s) -\phi_{1}^K(s)
\\
-(\phi_{1}^K)'(s) -(\phi_{1}^K)'(s)
\end{bmatrix},
\end{align*}
and
\begin{align*}
&\begin{bmatrix}
\cosh (n+1)s & \cosh(n-1)s
\\
(n+1)\sinh (n+1)s & (n-1)\sinh (n-1)s	
\end{bmatrix}
\begin{bmatrix}
c_n
\\
d_n	
\end{bmatrix}
=
\begin{bmatrix}
-\phi_{n}^K(s) -\phi_{n}^K(s)
\\
-(\phi_{n}^K)'(s) -(\phi_{n}^K)'(s)
\end{bmatrix},
\end{align*}
for $n\geq2$. Solving these linear systems yields the expressions given in the lemma for $c_n$ and $d_n$ in terms of $K_v$.
We then determine the constant $K_v$ by imposing the condition
\beq\label{eq:cndnsum_pre_cond}
c_0 + \sum_{n=1}^\infty (c_n + d_n) = 0.
\eeq
This condition is required to prove
\beq\label{eq:v2q2_decay}
(\Bv_2-\BU,q_2) \in \Mcal.
\eeq
We will be brief in presenting the proof of \eqnref{eq:v2q2_decay} since it is parallel to \eqnref{voneMcal}. We only mention why the condition \eqnref{eq:cndnsum_pre_cond} is required, and write down the formula for the pressure term $q_2$ since it will be used in latter part of this section.

Similarly to \eqnref{eq:anbnsum_conv}, one can show that for any positive number $k$ there is a constant $C$ such that
\beq\label{cndn}
\sum_{n=1}^\infty n^k e^{ns}(|c_n| + |d_n|) \le C.
\eeq
Note that
$$
(h \Phi^{e,F})(\Gz,\Gt) =c_0 \cosh \Gz + d_0 \Gz \sinh \Gz
+\sum_{n=1}^\infty \big( c_n \tilde{w}_n^+(\Gz,\Gt) +d_n \tilde{w}_n^-(\Gz,\Gt) \big) ,
$$
where
$$
\tilde{w}_n^\pm (\Gz,\Gt)= \cosh (n\pm 1)\Gz \cos n\Gt.
$$
Thanks to \eqnref{eq:cndnsum_pre_cond}, we have
$$
(h \Phi^{e,F})(\Gz,\Gt) =c_0 (\cosh \Gz-1) + d_0 \Gz \sinh \Gz
+\sum_{n=1}^\infty c_n (\tilde{w}_n^+(\Gz,\Gt)-1) +d_n (\tilde{w}_n^-(\Gz,\Gt)-1).
$$
We then use \eqnref{eq:wtn_estim} to obtain
$$
(h \Phi^{e,F})(\Gz,\Gt) = O(\Gz^2+\Gt^2).
$$

We use \eqref{eq:Laplacian_Psi} to see that $\GD\Phi_2^K=0$ and
\begin{align*}
    \GD \Phi_2 = \GD \Phi_2^F &= \frac{2}{a} c_0 + d_0 \frac{2}{a} (1-\tilde{w}_1)
     + \frac{2}{a}\sum_{n=1}^\infty((n+1)c_n - (n-1)d_n )\tilde{w}_n \\
     &\quad - \frac{2}{a}\sum_{n=1}^\infty n c_n \tilde{w}_{n+1}
     + \frac{2}{a}\sum_{n=1}^\infty n d_n \tilde{w}_{n-1}.
\end{align*}
Since the pressure $q_2$  is a harmonic conjugate of $\mu\GD\Phi^e$, we have
\begin{align}
    q_2 &=\frac{2\mu }{a} d_0 {w}_1 + \frac{2\mu}{a}\sum_{n=1}^\infty((n+1)c_n - (n-1)d_n ){w}_n \nonumber
    \\&\quad - \frac{2\mu}{a}\sum_{n=1}^\infty n c_n {w}_{n+1}
     + \frac{2\mu}{a}\sum_{n=1}^\infty n d_n{w}_{n-1} + C \label{qtwo}
\end{align}
for some constant $C$. We choose $C=0$. Then, since
\beq\label{wnest}
|w_n(\Gz,\Gt)| \lesssim n^2 e^{n s}(\Gz^2 + \Gt^2),
\eeq
we have $q_2 = O(\Gz^2 + \Gt^2)$ as $(\Gz, \Gt) \to 0$, namely, \eqnref{eq:v2q2_decay} holds.
\qed

\subsubsection{Stream function for $(\Bh_{\mathrm{rot}},p_{\mathrm{rot}})$}

\begin{lemma}\label{lem:Psi_rot_exp}
The stream function $\Phi_{\mathrm{rot}}$ associated with the solution $(\Bh_{\mathrm{rot}},p_{\mathrm{rot}})$ is given by
\begin{align}
(h\Phi_{\mathrm{rot}})(\Gz,\Gt)
&= K_{\mathrm{rot}}(\cosh\Gz -\cos\Gt)\ln (2\cosh \Gz - 2\cos \Gt) +a_0' \cosh \Gz + d_0' \Gz \sinh \Gz \nonumber
\\
&\quad
+\sum_{n=1}^\infty \big(a_n' \cosh (n+1) \Gz + b_n' \cosh (n-1)\Gz \big)\cos n\Gt ,
\label{eq:Psi_rot}
\end{align}
where
\begin{align*}
K_{\mathrm{rot}}=-a  \left(\frac{s \sinh^2 s \tanh s}{\sinh s \cosh s +s}+\frac{1}{2} + M\right)^{-1}
\end{align*}
with $M$ given in \eqnref{M}, and
\begin{align*}
    a_0' &= a-\frac{   K_{\mathrm{rot}}(s^2+s+e^{-s}\sinh s)}{\sinh s \cosh s +s},
    \quad
    d_0' =  -\frac{ K_{\mathrm{rot}} \sinh^2 s}{\sinh s \cosh s +s} ,
    \\
    a_1' &= \frac{1}{2} K_{\mathrm{rot}} e^{-s} \mathrm{sech}\, s,
    \quad
    b_1' = K_{\mathrm{rot}}(s+1-\frac{1}{2} \tanh s),
    \\
    a_n' &=\frac{2K_{\mathrm{rot}} (n e^{-s} \sinh s +e^{-n s}\sinh n s)}{n(n+1)(\sinh 2ns + n \sinh 2s)},
    \\
    b_n' &=-\frac{2K_{\mathrm{rot}} (n e^{s} \sinh s +e^{-n s}\sinh n s)}{n(n-1)(\sinh 2ns + n \sinh 2s)}.
\end{align*}
\end{lemma}

\proof
Let
$$
(\tilde\Bh_\mathrm{rot},\tilde{p}_\mathrm{rot}):=(\Bh_{\mathrm{rot}},p_{\mathrm{rot}}) - (\Bpsi_3,0).
$$
Then $(\tilde\Bh_\mathrm{rot},\tilde{p}_\mathrm{rot})$ is the solution to
\beq
\begin{cases}
\mu \GD \tilde\Bh_\mathrm{rot} = \nabla \tilde{p}_\mathrm{rot} \quad &\mbox{in }D^e,
\\
\nabla \cdot \tilde\Bh_\mathrm{rot} = 0 \quad &\mbox{in }D^e,
\\
(\tilde\Bh_\mathrm{rot},\tilde{p}_\mathrm{rot}) - (-\Bpsi_3,0) \in\mathcal{M},
\end{cases}
\eeq
with the no-slip boundary condition, namely,
\beq\label{eq:stream_htilde_rot_bcbc}
\tilde\Bh_\mathrm{rot} |_{\p D_1}=0, \quad \tilde\Bh_\mathrm{rot}|_{\p D_2}=0.
\eeq
Observe that the above equation is similar to the equation \eqref{Bv2} for $(\Bv_2,q_2)$ with the only difference being that the background solution $\BUsh$ is replaced with $-\Bpsi_3$.

It is easy to see that the function $\Phi_\mathrm{rot}^0$ defined by
$$
\Phi_\mathrm{rot}^0 = -\frac{1}{2}(x^2+y^2)
$$
is a stream function associated with the solution $(-\Bpsi_3,0)$. In bipolar coordinates,
\beq
\Phi_\mathrm{rot}^0 = \frac{1}{2}\frac{a^2 (-\sinh^2 \Gz+ \sin^2 \Gt)}{(\cosh\Gz-\cos\Gt)^2}.
\eeq
Note that $\Phi_\mathrm{rot}^0$ has the even symmetry in both $\Gz$ and $\Gt$.
In exactly the same way as in the proof of Lemma \ref{lem:Phi_2}, we can find the stream function associated with $(\tilde\Bh_\mathrm{rot},\tilde{p}_\mathrm{rot})$, which immediately yields Lemma \ref{lem:Psi_rot_exp}.
\qed

\subsection{Asymptotics of $K_v$ and $K_{\mathrm{rot}}$}

\begin{lemma}\label{lem:K_asymp}
As $\Gd\rightarrow 0$, we have
\begin{align}
K_{v} &=R\frac{1-G_0}{F_0} \sqrt{\frac{R}{\Gd}} + O(1), \label{eq:Kv_asymp} \\
K_{\mathrm{rot}} &= -\frac{R}{F_0} \sqrt{\frac{R}{\Gd}}+O(1),\label{eq:Krot_asymp}
\end{align}
where $F_0$ and $G_0$ are the numbers given in \eqref{F0G0}.
\end{lemma}
\proof
The proof is based on a special case of the Euler-Maclaurin summation formula: if $f\in C^1(\Rbb^+) \cap L^1(\Rbb^+)$, then, for a small parameter $s>0$, we have
\beq\label{EM}
s \sum_{n=0}^\infty f(x_0+ n s) = \int_{x_0}^\infty f(x) dx + R_1,
\eeq
where the remainder term $R_1$ satisfies
$$
|R_1| \lesssim s \left( |f(x_0)| + \int_{x_0}^\infty |f'(x)|dx \right).
$$

We first consider the asymptotics of the series $M$ defined by (\ref{M}).
One can easily see that
\begin{align*}
    M + 2\sum_{n=2}^\infty \frac{1}{n(n^2-1)} =  \sum_{n=2}^\infty \frac{4e^{-ns} \sinh^2(ns)(\cosh ns + \sinh ns) - 4n^2 \sinh^2 s }{n(n^2-1)(\sinh 2n s + n \sinh 2s)},
\end{align*}
and
$$
2\sum_{n=2}^\infty \frac{1}{n(n^2-1)} = \frac{1}{2}.
$$
Thus, the Euler-Maclaurin summation formula yields
\begin{align}
M + \frac{1}{2} &=  s^3 \sum_{n=2}^\infty \frac{ 4e^{-ns} \sinh^2(ns)(\cosh ns + \sinh ns) - 4(n s)^2 (\sinh s/s)^2  }{(ns)((ns)^2-s^2)(\sinh 2n s + ns (\sinh 2s/s)} \nonumber
\\
&=s^2\int_{2s}^\infty f_s(x) dx + s^2R_1,
\end{align}
where
$$
f_s(x) : = \frac{ 4e^{-x} \sinh^2 x (\cosh x + \sinh x)-4x^2(\sinh s /s)^2 }{x(x^2-s^2)(\sinh 2x + (\sinh 2s/s)x)},
$$
and
$$
|R_1| \lesssim s \left( |f_s(2s)| + \int_{2s}^\infty |f_s'(x)|dx \right).
$$
By straightforward but tedious computations, one can see that
$$
|f_s(2s)|\leq C, \quad \int_{2s}^\infty|f_s'(x)|dx \leq C',
$$
where
$C$ and $C'$ are constants independent of $s>0$.
Therefore, as $s\rightarrow 0$, we obtain
$$
M + \frac{1}{2} = s^2 \int_0^\infty f_0(x) dx  + O(s^3) =s^2F_0 + O(s^3).
$$
So, for small $s$, we have
$$
K_{\mathrm{rot}} = -a \left(  \frac{1}{2} + M + O(s^3) \right)^{-1} =  -\frac{a}{s^2 F_0} + O(1).
$$
The other quantity $K_v$ can be estimated similarly, and the proof is completed.
\qed

\section{No blow-up with no-slip boundary conditions II}\label{sec:noslip2}

\subsection{Proof of Theorem \ref{thm:boundedstress_v1q1}}

We first estimate the strain tensor $\Ecal[\Bv_1]$. Since $\Ecal[\BU]=O(1)$, we estimate $\Ecal[\Bv_1-\BU]$.

The following formulae are derived using the relations \eqref{eq:strain_bipolar1}--\eqref{eq:strain_bipolar3} between the strain tensor and the stream function $\Phi_1$ and \eqnref{eq:stream_v1}:
\begin{align*}
    \Ecal_{\Gz\Gz}[\Bv_1-\BU] &=
    -h(\Gz,\Gt) 2 a_1 \cosh 2\Gz \cos \Gt - h(\Gz,\Gt) b_1 \cos\Gt
    \\
    &\quad - h(\Gz,\Gt)\sum_{n=2}^\infty \Big(\tilde{a}_n \cosh (n+1)\Gz +\tilde{b}_n \cosh (n-1)\Gz\Big)\cos n\Gt, \\
    \Ecal_{\Gz\Gt}[\Bv_1-\BU] &= h(\Gz,\Gt) 2a_1 \sinh 2\Gz \sin \Gt
    \\
    &\quad + h(\Gz,\Gt)\sum_{n=2}^\infty \Big(\tilde{a}_n \sinh (n+1)\Gz +\tilde{b}_n \sinh (n-1)\Gz\Big)\sin n\Gt, \\
    \Ecal_{\Gt\Gt}[\Bv_1-\BU] &= - \Ecal_{\Gz\Gz}[\Bv_1-\BU] ,
\end{align*}
where
\beq\label{tildean}
\tilde{a}_n = n(n+1) a_n, \quad \tilde{b}_n = n(n-1) b_n.
\eeq
Here, $a_n$ and $b_n$ are given in Lemma \ref{lem:Phi_1}.

Using the hyperbolic identities
\begin{align*}
   & \cosh(n+1)\Gz + \cosh(n-1)\Gz = 2\cosh n \Gz \cosh \Gz,
\\
    &\cosh(n+1)\Gz - \cosh(n-1)\Gz = 2\sinh n \Gz \sinh \Gz,
\end{align*}
we can rewrite $\Ecal_{\Gz\Gz}$ and $ \Ecal_{\Gz\Gt}$ as
\begin{align}
    & \Ecal_{\Gz\Gz}[\Bv_1-\BU] =
     -h(\Gz,\Gt) 2 a_1 \cosh 2\Gz \cos \Gt - h(\Gz,\Gt) b_1 \cos\Gt\nonumber
    \\
    &\quad - h(\Gz,\Gt)\sum_{n=2}^\infty \big( ({\tilde{a}_n+\tilde{b}_n}) \cosh n \Gz \cosh \Gz + ({\tilde{a}_n- \tilde{b}_n}) \sinh n\Gz \sinh \Gz \big)\cos n\Gt,
    \label{eq:Ecal_GzGz_v1_symm}
\end{align}
and
\begin{align}
    &\Ecal_{\Gz\Gt}[\Bv_1-\BU] =
     h(\Gz,\Gt) 2a_1 \sinh 2\Gz \sin \Gt \nonumber
    \\
    &\quad + h(\Gz,\Gt)\sum_{n=2}^\infty \big( ({\tilde{a}_n+\tilde{b}_n}) \sinh n \Gz \cosh \Gz + ({\tilde{a}_n- \tilde{b}_n}) \cosh n\Gz \sinh \Gz \big)\sin n\Gt.
    \label{eq:Ecal_GzGt_v1_symm}
\end{align}

From the expressions of $a_n$ and $b_n$ given in Lemma \ref{lem:Phi_1}, we have, for $n \geq 2$,
\begin{align*}
{\tilde{a}_n+\tilde{b}_n} &=-\frac{4a}{s}\frac{n s e^{-n s} \sinh n s -   (n s)^2\eta_2 +   (n s)^3 \eta_1}{\sinh 2ns - 2ns\eta_{2}},
\\
{\tilde{a}_n-\tilde{b}_n} &=-\frac{4a}{s^2}\frac{(ns)^2 e^{-n s} \sinh n s -   (n s)^3 \eta_2 + s^2   (n s)^2 \eta_1}{\sinh 2ns - 2ns\eta_{2}},
\end{align*}
where
\beq\label{eq:def_eta1_eta2}
\eta_1 = \eta_{1}(s) := \frac{\sinh^2 s}{s^2} , \quad \eta_2 = \eta_{2}(s) := \frac{\sinh 2s}{2s}.
\eeq
If we define
\begin{align}
    f_1(x) := \frac{x e^{-x}\sinh x  - x^2\eta_{2} + x^3 \eta_1}{\sinh 2x - 2x\eta_{2}}, \quad
f_2(x) := \frac{x^2 e^{-x}\sinh x  -  x^3 \eta_2 + s^2 x^2  \eta_1 }{\sinh 2x - 2x\eta_{2}},\label{eq:def_f1f2}
\end{align}
then ${\tilde{a}_n+\tilde{b}_n}$ and ${\tilde{a}_n-\tilde{b}_n}$ can be rewritten as
\begin{align}
{\tilde{a}_n+\tilde{b}_n} =-\frac{4a}{s} f_1(ns), \quad
{\tilde{a}_n-\tilde{b}_n} =-\frac{4a}{s^2} f_2(ns) .\label{eq:anpmbn_fj}
\end{align}

It follows from $\Ecal_{\Gz\Gz}$ and $\Ecal_{\Gz\Gt}$ that
from \eqref{eq:Ecal_GzGz_v1_symm} and \eqref{eq:anpmbn_fj} that
\begin{align}
     \Ecal_{\Gz\Gz}[\Bv_1-\BU] &=
     -h(\Gz,\Gt) 2 c_1 \cosh 2\Gz \cos \Gt - h(\Gz,\Gt) d_1 \cos\Gt\nonumber
    \\
    &\quad + \frac{4a}{s}h(\Gz,\Gt)\sum_{n=2}^\infty \big(  f_1(ns) \cosh n \Gz \cosh \Gz \big)\cos n\Gt\nonumber
    \\
    &\quad + \frac{4a}{s^2}h(\Gz,\Gt)\sum_{n=2}^\infty \big(   f_2(ns) \sinh n\Gz \sinh \Gz \big)\cos n\Gt. \label{eq:Ecal_GzGz_pre_series}
\end{align}
Let, for $j=1,2$,
\beq
A_{j,n}^+(\Gz)  = f_j(n s) \cosh n\Gz, \quad A_{j,n}^-(\Gz)  = f_j(n s) \sinh n\Gz,
\eeq
and then define $S_j^{++}$, $S_j^{-+}$, etc, by
\beq
S_j^{\pm+} = \sum_{n=2}^\infty A_{j,n}^\pm(\Gz)  ah(\Gz,\Gt)\cos n\Gt, \quad
S_j^{\pm-} = \sum_{n=2}^\infty A_{j,n}^\pm(\Gz)  ah(\Gz,\Gt)\sin n\Gt.
\label{eq:Sj_pm_rewrite}
\eeq
Then, \eqref{eq:Ecal_GzGz_pre_series} reads
\begin{align}
    \Ecal_{\Gz\Gz}[\Bv_1-\BU]
    &=S_0 + \frac{4}{s} \cosh \Gz S_1^{++} + \frac{4}{s^2} \sinh \Gz S_2^{-+},
    \label{eq:Ecal_GzGz_v1_preasymp}
\end{align}
where
\beq\label{eq:def_S0}
S_0 = -h(\Gz,\Gt) 2 a_1 \cosh 2\Gz \cos \Gt - h(\Gz,\Gt) b_1 \cos\Gt.
\eeq
Similarly, one can see that $\Ecal_{\Gz\Gt}$ is written as
\begin{align}
    \Ecal_{\Gz\Gt}[\Bv_1-\BU]
    &=\widetilde{S}_0 - \frac{4}{s} \cosh \Gz S_1^{--} - \frac{4}{s^2} \sinh \Gz S_2^{+-},
    \label{eq:Ecal_GzGt_v1_preasymp}
\end{align}
where
\beq\label{eq:def_S0p}
\widetilde{S}_0 =  h(\Gz,\Gt) 2b_1 \sinh 2\Gz \sin \Gt.
\eeq

We use the following lemma here and present its proof in Appendix \ref{appendixD}.
\begin{lemma}\label{lem:fj_v1_estim}
If $2s \le x$, then
\begin{equation}\label{eq:fj_estim}
|f_j(x)|+ |f_j'(x)| + |f_j''(x)|\lesssim (1+x^3)e^{-2x}, \quad j=1,2.
\end{equation}
If $2s \le x \le 1$, then
\beq\label{eq:Df1_DDf2_asymp}
\left| f_1'(x) - \frac{1}{2} \right| \lesssim x,
\quad
|f_2''(x) - {1}| \lesssim x.
\eeq
\end{lemma}

\begin{lemma}\label{lem:pre_asymp_Sj}
The following asymptotic formulas hold for $j=1,2$:
\begin{align}
S_{j}^{+ +} &=- \frac{ 1}{2} f_j(2s)\cosh 2\Gz \cos \Gt  + f_j(2s)\cosh 2\Gz \cos 2\Gt \nonumber \\
&\quad -\frac{1}{2} f_j(3s)\cosh 3\Gz \cos 2\Gt  + O(s), \label{Sj++}\\
S_{j}^{+ -}&=- \frac{ 1}{2} f_j(2s)\cosh 2\Gz \sin \Gt  +f_j(2s)\cosh 2\Gz \sin 2\Gt \nonumber \\
&\quad  -\frac{1}{2} f_j(3s)\cosh 3\Gz \sin 2\Gt  + O(s), \label{Sj+-} \\
S_{j}^{- \pm} &=O(s), \label{Sj-pm}
\end{align}
as $s \to 0$.
\end{lemma}

\proof
We first have the following identity:
\begin{align*}
ah(\Gz,\Gt)\cos n\Gt &= \cosh\Gz \cos n \Gt - \cos\Gt \cos n\Gt
\\
&= -\frac{1}{2}\big(\cos (n+1)\Gt - 2\cosh \Gz \cos n \Gt + \cos (n-1)\Gt \big)
\\
&=-\frac{1}{2}\big(\cos (n+1)\Gt - 2 \cos n \Gt + \cos (n-1)\Gt \big) + (\cosh\Gz-1) \cos n \Gt
\\
&=-\frac{1}{2}\big(\cos (n+1)\Gt - 2 \cos n \Gt + \cos (n-1)\Gt \big) + \sinh^2(\Gz/2) \cos n \Gt.
\end{align*}
By substituting this identity into \eqref{eq:Sj_pm_rewrite} and then rearranging indices, we arrive at
\begin{align}
 S_j^{\pm+} &=
 - \frac{ 1}{2} (A_{j,2}^\pm \cos \Gt  -2 A_{j,2}^\pm \cos 2\Gt + A_{j,3}^\pm \sin 2\Gt )
 \nonumber \\ &\quad
 -\frac{1}{2}\sum_{n=3}^\infty  \Big(A_{j,n+1}^\pm - 2A_{j,n}^\pm + A_{j,n-1}^\pm   \Big)\cos n\Gt \nonumber
 \\
 & \quad +  \sinh^2 (\Gz/2) \sum_{n=2}^\infty A_{j,n}^\pm \cos n\Gt
=: S_{j,0}^{\pm +} + S_{j,1}^{\pm +} + S_{j,2}^{\pm +}. \label{eq:v2_Sp_fourier}
\end{align}

Note that $A_{j,n}^\pm$ is of the form
$$
A_{j,n}^\pm = F_{j,n} G^\pm_{n},
$$
where
$$
F_{j,n} = f_j(n s), \quad G_n^+ = \cosh n \Gz, \quad G_n^- = \sinh n\Gz.
$$
We then have
\begin{align*}
     A_{j,n+1}^\pm  - 2A_{j,n}^\pm + A_{j,n-1}^\pm  & = F_{j,n+1} G_{n+1}^\pm - 2 F_{j,n} G_n^\pm + F_{j,n-1}G_{n-1}^\pm
    \\
     &= (F_{j,n+1}-2F_{j,n}+F_{j,n-1}) G_n^\pm  + F_{j,n} (G^\pm_{n+1}-2G^\pm_n + G^\pm_{n-1})
    \\
    &+ (F_{j,n+1} - F_{j,n}) (G^\pm_{n+1}- G^\pm_{n})
     + (F_{j,n} - F_{j,n-1}) (G^\pm_{n}- G^\pm_{n-1}).
\end{align*}
One can easily see that
\begin{align*}
&G_{n+1}^\pm - G_n^\pm = \sinh ({\Gz}/{2}) (e^{(n+\frac{1}{2})\Gz} \mp e^{-(n+\frac{1}{2})\Gz} ),
\\
    &G_{n+1}^\pm -2G_n^\pm + G_{n-1}^\pm =  2\sinh^2 (\Gz/2) (e^{n\Gz} \pm e^{-n\Gz} ).
\end{align*}
Since $|\Gz| \le s$, we have
\beq\label{Gnest}
|G_n^\pm| \lesssim e^{ns},\quad |G_{n+1}^\pm - G_n^\pm| \lesssim s e^{n s},
\quad
 | G_{n+1}^\pm -2G_n^\pm + G_{n-1}^\pm | \lesssim s^2 e^{n s}.
\eeq

Next, we estimate $F_{j,n}$ and its finite differences.
By the mean value theorem, there exist $x_n^*\in(ns, (n+1)s), x_n^{**}\in ((n-1)s, (n+1)s)  $ such that
$$
\frac{F_{j,n+1}-F_{j,n}}{s} = f_j'(x_n^*), \quad \frac{F_{j,n+1}-2F_{j,n}+F_{j,n-1}}{s^2} = f_j''( x_n^{**} ).
$$
Then, by \eqnref{eq:fj_estim}, we infer
\begin{align}
&|F_{j,n} | \lesssim ( 1+ (n s)^3) e^{-2ns},
\\
&|F_{j,n+1}-F_{j,n}| \lesssim  s ( 1+ (n s)^3) e^{-2ns},
\\
&|F_{j,n+1}-2F_{j,n} + F_{j,n-1}| \lesssim  s^2 ( 1+ (n s)^3) e^{-2ns}.
\label{eq:Fn_estim}
\end{align}

These estimates together with \eqnref{Gnest} lead us to
$$
|A_n^\pm |\lesssim e^{n s}|F_{j,n}| \lesssim (1+(ns)^3) e^{- n s},
$$
and
\begin{align*}
     |A_{n+1}^\pm  - 2A_{n}^\pm + A_{n-1}^\pm|  & \lesssim
      \Big|F_{n+1}-2F_n+F_{n-1}\Big| e^{n s}  + F_n s^2 e^{n s}
    \\
    &\quad + |F_{n+1} - F_{n}| s e^{n s}  + |F_n - F_{n-1}| s e^{ns}
    \\
    &\lesssim s^2 (1+(ns)^3)e^{-n s}.
\end{align*}
Using \eqnref{EM}, we have
\begin{align*}
    |S_{j,1}^{\pm+}| &\lesssim \sum_{n=2} |A_{n+1}^\pm  - 2A_{n}^\pm + A_{n-1}^\pm| \lesssim \sum_{n=2} s^2 (1+(ns)^3)e^{-ns}
    \\
    &\lesssim s\int_0^\infty (1+x^3)e^{-x} dx \lesssim s,
\end{align*}
and
\begin{align*}
    |S_{j,2}^{\pm+}| &\lesssim s^2 \sum_{n=2} | A_{n}^\pm|\lesssim s^2\sum_{n=2}  (1+(ns)^3)e^{-ns}
    \\
    &\lesssim s\int_0^\infty (1+x^3)e^{-x} dx \lesssim s.
\end{align*}
Therefore, from \eqref{eq:v2_Sp_fourier} and \eqref{eq:v2_Sm_fourier}, we see that
\begin{align*}
S_{j}^{+ +}&= S_{j,0}^{+ +} + O(s),
\end{align*}
which is the formula \eqnref{Sj++}.
Similarly,
\begin{align*}
S_{j}^{- +}&= S_{j,0}^{- +} + O(s)
\\
&=- \frac{ 1}{2} (A_{j,2}^- \cos \Gt  -2 A_{j,2}^- \cos 2\Gt + A_{j,3}^- \cos 2\Gt ) + O(s)
\\
&=- \frac{ 1}{2} (f_j(2s)\sinh 2\Gz \cos \Gt  -2 f_j(2s)\sinh 2\Gz \cos 2\Gt + f_j(3s)\sinh 3\Gz \cos 2\Gt ) + O(s).
\end{align*}
Since $\sinh \Gz=O(s)$ and $|f_j(2s)|+|f_j(3s)|$ is bounded thanks to \eqnref{eq:fj_estim}, the estimate \eqnref{Sj-pm} for $S_j^{-+}$ follows.

Using the identity
\begin{align*}
	ah(\Gz,\Gt)\sin n\Gt =  -\frac{1}{2}\big(\sin (n+1)\Gt - 2 \sin n \Gt + \sin (n-1)\Gt \big) + \sinh^2(\Gz/2) \sin n \Gt,
\end{align*}
one can see that
\begin{align}
 S_j^{\pm-} &=
 - \frac{ 1}{2} (A_{j,2}^\pm \sin \Gt  -2 A_{j,2}^\pm \sin 2\Gt + A_{j,3}^\pm \sin 2\Gt )
 \nonumber \\ &\quad
 -\frac{1}{2}\sum_{n=3}^\infty  \Big(A_{j,n+1}^\pm - 2A_{j,n}^\pm + A_{j,n-1}^\pm   \Big)\sin n\Gt +  \sinh^2 (\Gz/2) \sum_{n=2}^\infty A_{j,n}^\pm \sin n\Gt.
    \label{eq:v2_Sm_fourier}
\end{align}
The other formulas, namely, \eqnref{Sj+-} and \eqnref{Sj-pm} for $S_j^{--}$, can be proved in the same way.
The proof is completed.
\qed

\medskip

We are now prepared for estimating $\Ecal_{\Gz\Gz}$ and $\Ecal_{\Gz\Gt}$. By applying Lemma \ref{lem:pre_asymp_Sj} to \eqref{eq:Ecal_GzGz_v1_preasymp} and \eqref{eq:Ecal_GzGt_v1_preasymp}, we have
\begin{align*}
    \Ecal_{\Gz\Gz}[\Bv_1-\BU]
    &=S_0 + \frac{4}{s} \cosh \Gz S_1^{++} + O(1)
    \\
    &=S_0 - \frac{2}{s} \cosh \Gz \Big[f_1(2s)\cosh 2\Gz \cos \Gt   \\
    &\quad -2 f_1(2s)\cosh 2\Gz \cos 2\Gt+ f_1(3s)\cosh 3\Gz \cos 2\Gt \Big] +O(1),
  \end{align*}
and
\begin{align*}
    \Ecal_{\Gz\Gt}[\Bv_1-\BU]
    &=\widetilde{S}_0 - \frac{4}{s^2} \sinh \Gz S_2^{+-} + O(1)
    \\
    &= \widetilde{S}_0 +\frac{2}{s^2}\sinh \Gz \Big[f_j(2s)\cosh 2\Gz \sin \Gt
    \\
    &\quad -2 f_j(2s)\cosh 2\Gz \sin 2\Gt + f_j(3s)\cosh 3\Gz \sin 2\Gt \Big] +O(1).
\end{align*}
By applying Taylor expansions to $f_j$ given in \eqref{eq:def_f1f2}, we see that
\begin{align*}
&f_1(2s) = s + O(s^2), \quad f_1(3s) = \frac{3}{2} s + O(s^2),
\\
&f_2(2s) = -\frac{3}{2}s + O(s^2), \quad f_2(3s) = -\frac{9}{4} s + O(s^2).
\end{align*}
So we have
\begin{align*}
	\Ecal_{\Gz\Gz}[\Bv_1] &=S_0 + O(1), \qquad \Ecal_{\Gz\Gt}[\Bv_1] =\widetilde{S}_0 +O(1).
\end{align*}

It remains to estimate $S_0$ and $\widetilde{S}_0$.
 Applying Taylor expansions to $a_1$ and $b_1$ given in Lemma \ref{lem:Phi_1}, we have
\beq\label{a1b1}
 a_1  = \frac{3a}{4s} +O(s), \quad b_1  = -\frac{3a}{2s} +O(s).
\eeq
Then, from \eqref{eq:def_S0} and \eqref{eq:def_S0p}, we have
\begin{align*}
S_0 &= -h(\Gz,\Gt) (2 a_1 \cosh 2\Gz + b_1 )\cos\Gt = - h(\Gz,\Gt)(2a_1 +b_1) + O(1) = O(1),
\\
\widetilde{S}_0 &=  h(\Gz,\Gt) 2a_1 \sinh 2\Gz \sin \Gt = O(1).	
\end{align*}
Therefore, we obtain that
\begin{align}
		\Ecal_{\Gz\Gz}[\Bv_1] &= O(1), \qquad \Ecal_{\Gz\Gt}[\Bv_1] = O(1).
\end{align}

We now prove that the pressure $q_1$ is bounded regardless of $\Gd$. Recall from \eqnref{qone2} that $q_1$ is given by
\begin{align*}
q_1 &=-a_1\frac{2\mu}{a} (2 (\tilde{w}_1-1)-(\tilde{w}_2-1) ) +b_1\frac{2\mu}{a}(\tilde{w}_1-1)\nonumber
 \\
 &\quad - \frac{2\mu}{a}\sum_{n=2}^\infty((n+1)a_n - (n-1)b_n )(\tilde{w}_n -1)
 \nonumber   \\&\quad
    +     \frac{2\mu}{a}\sum_{n=2}^\infty n a_n (\tilde{w}_{n+1}-1)
     - \frac{2\mu}{a}\sum_{n=2}^\infty n b_n (\tilde{w}_{n-1}-1),
\end{align*}
where $\tilde{w}_n(\Gz,\Gt)=\cosh n\Gz \cos n \Gt$. Using notation in \eqnref{tildean}, we have
\begin{align}
q_1 &=\frac{2\mu}{a} (-2a_1 + b_1-2 b_2) \tilde{w}_1 + \frac{2\mu}{a} (a_1 -3 a_2 + b_2 - 3 b_3) \tilde{w}_2 \nonumber
        \\
&\quad -\frac{2\mu}{a}\sum_{n=3}^\infty \frac{1}{n}(\tilde{a}_n - \tilde{b}_n - \tilde{a}_{n-1} + \tilde{b}_{n+1}) \tilde{w}_n.
    \label{q2_series}
\end{align}
Note that
\begin{align*}
    \tilde{a}_n - \tilde{b}_n - \tilde{a}_{n-1} + \tilde{b}_{n+1} &=\frac{1}{2} \Big((\tilde{a}_{n+1}+\tilde{b}_{n+1})  -(\tilde{a}_{n-1}+\tilde{b}_{n-1}) \Big)
    \\
    &\quad - \frac{1}{2}\Big((\tilde{a}_{n+1}-\tilde{b}_{n+1})-2(\tilde{a}_n - \tilde{b}_n)   + (\tilde{a}_{n-1}-\tilde{b}_{n-1}) \Big).
\end{align*}
Then we have from \eqref{eq:anpmbn_fj} that
\begin{align*}
    \tilde{a}_n - \tilde{b}_n - \tilde{a}_{n-1} + \tilde{b}_{n+1} &= -\frac{2a}{s} \big( f_1((n+1)s)- f_1((n-1)s)\big)
    \\
    &\quad   +\frac{2a}{s^2} \big(f_2((n+1)s)-2f_2(ns) + f_2((n-1)s)\big).
\end{align*}
Therefore, by \eqref{q2_series}, we have
\begin{align}
|q_1| &\lesssim \frac{1}{s} |-2a_1 + b_1-2 b_2| + \frac{1}{s} |a_1 -3 a_2 + b_2 - 3 b_3|\nonumber
        \\
        &\quad +\sum_{n=3}^\infty \frac{s\cosh ns}{ns}\Big|2\frac{f_1((n+1)s)- f_1((n-1)s)}{2s} \nonumber
        \\
        &\qquad\quad\qquad -
  \frac{1}{s^2} \big(f_2((n+1)s)-2f_2(ns) + f_2((n-1)s)\big)\Big|. \label{eq:q1_prefinal_estim}
\end{align}
By applying the mean value theorem, we have
\begin{align}
|q_1| &\lesssim 	 \frac{1}{s} |-2a_1 + b_1-2 b_2| + \frac{1}{s} |a_1 -3 a_2 + b_2 - 3 b_3|\nonumber
        \\
        &\quad +\sum_{n=3}^\infty \frac{s\cosh ns}{ns}\Big|2 f_1'(x_n^*) -f_2''(x_n^{**})\Big| \label{eq:q1_pre_estim}
\end{align}
for some $x_n^* \in ((n-1)s,(n+1)s)$ and $x_n^{**}\in ((n-1)s,(n+1)s)$.

By regarding the infinite series in \eqnref{eq:q1_pre_estim} as a Riemann sum, we infer
\begin{align*}
I:= \sum_{n=3}^\infty \frac{s\cosh ns}{ns}\Big|2 f_1'(x_n^*) -f_2''(x_n^{**})\Big|
 \lesssim \int_0^\infty \frac{\cosh x}{x} |2f_1'(x) - f_2''(x)| dx.
\end{align*}
It then follows that
\begin{align*}
I & \lesssim \int_0^1 + \int_1^\infty \frac{\cosh x}{x} |2f_1'(x) -f_2''(x)| dx \\
& = \int_0^1 \frac{\cosh x}{x} |2(f_1'(x)-\frac{1}{2}) - (f_2''(x)-1)| dx + \int_1^\infty \frac{\cosh x}{x} |2f_1'(x) -f_2''(x)| dx \\
& \lesssim \int_0^1 \frac{\cosh x}{x} (x+x) dx + \int_1^\infty \frac{\cosh x}{x} x^3 e^{-2x} dx \\
	& \lesssim 1 + \int_1^\infty x^2 e^{-x} dx \lesssim 1,
\end{align*}
where the third inequality follows from \eqnref{eq:Df1_DDf2_asymp}.

We next estimate the first two terms in the right-hand side of \eqref{eq:q1_prefinal_estim}. By Taylor expansions we obtain
\beq\label{eq:a1a2b1b2b3_estim}
a_2  = \frac{1}{2} \frac{a}{s} + O(s), \quad b_2  = -\frac{3}{2} \frac{a}{s} + O(s),
\quad b_3  = -\frac{3}{4} \frac{a}{s} + O(s).
\eeq
These together with \eqnref{a1b1} yield
\begin{align*}
-2 a_1 + b_1 -2 b_2 = O(s), \quad a_1-3a_2+ b_2 - 3 b_3 = O(s).
\end{align*}
Therefore, from \eqref{eq:q1_pre_estim}, we have
\begin{align*}
|q_1|\lesssim 1.	
\end{align*}
This completes the proof. \qed

\subsection{Proofs of Theorems \ref{thm:boundedstress_rot} and \ref{thm:boundedstress_v2q2}}

We only prove Theorem \ref{thm:boundedstress_v2q2}. Thanks to the similarity between the stream functions $\Phi_{2}$ and $\Phi_{\mathrm{rot}}$, Theorem \ref{thm:boundedstress_rot} can be proved in exactly the same way.

We first estimate the strain tensor $\Ecal[\Bv_2-\BU]$. In this proof, $\BU=\BUsh$.
Using \eqref{eq:strain_bipolar1}-\eqref{eq:strain_bipolar3} and \eqnref{eq:W_biharmonic_solution}, one can see that
\begin{align}
    \Ecal_{\Gz\Gz}[\Bv_2-\BU] &=
    -K_v \frac{\sinh\Gz}{a}\sin \Gt +  h(\Gz,\Gt)2 c_1 \sinh 2\Gz \sin \Gt\nonumber \\
    &\quad + h(\Gz,\Gt)\sum_{n=2}^\infty \Big(\tilde{c}_n \sinh (n+1)\Gz +\tilde{d}_n \sinh (n-1)\Gz\Big)\sin n\Gt, \\
    \Ecal_{\Gz\Gt}[\Bv_2-\BU] &= K_v \frac{\cosh2\Gz -2\cosh\Gz\cos\Gt+ \cos 2\Gt}{2a} \nonumber
    \\
    &\quad +  h(\Gz,\Gt) d_0\cosh\Gz + h(\Gz,\Gt) 2c_1\cosh 2\Gz \cos\Gt\nonumber
    \\
    &\quad + h(\Gz,\Gt)\sum_{n=2}^\infty \Big(\tilde{c}_n \cosh (n+1)\Gz +\tilde{d}_n \cosh (n-1)\Gz\Big)\cos n\Gt, \\
    \Ecal_{\Gt\Gt}[\Bv_2-\BU] &= -\Ecal_{\Gz\Gz}[\Bv_2-\BU],
\end{align}
where
$$
\tilde{c}_n = n(n+1) c_n, \quad \tilde{d}_n = n(n-1) b_n.
$$

Using the hyperbolic identities
\begin{align*}
    &\sinh(n+1)\Gz + \sinh(n-1)\Gz = 2\sinh n \Gz \cosh \Gz,
\\
    &\sinh(n+1)\Gz - \sinh(n-1)\Gz = 2\cosh n \Gz \sinh \Gz,
\end{align*}
we can rewrite $\Ecal_{\Gz\Gz}$ and $ \Ecal_{\Gz\Gt}$ as
\begin{align}
    &\Ecal_{\Gz\Gz}[\Bv_2-\BU] =
    -K_v \frac{\sinh\Gz}{a}\sin \Gt +  h(\Gz,\Gt)2 c_1 \sinh 2\Gz \sin \Gt\nonumber
    \\
    &\quad + h(\Gz,\Gt)\sum_{n=2}^\infty \big( ({\tilde{c}_n+\tilde{d}_n}) \sinh n \Gz \cosh \Gz + ({\tilde{c}_n- \tilde{d}_n}) \cosh n\Gz \sinh \Gz \big)\sin n\Gt, \label{eq:Ecal_GzGz_v2_symm}
\end{align}
and
\begin{align}
    &\Ecal_{\Gz\Gt}[\Bv_2-\BU] =
    K_v \frac{\cosh2\Gz -2\cosh\Gz\cos\Gt+ \cos 2\Gt}{2a}  \nonumber
    \\
    &\quad +  h(\Gz,\Gt) d_0\cosh\Gz + h(\Gz,\Gt) 2c_1\cosh 2\Gz \cos\Gt\nonumber
    \\
    &\quad + h(\Gz,\Gt)\sum_{n=2}^\infty \big( ({\tilde{c}_n+\tilde{d}_n}) \cosh n \Gz \cosh \Gz + ({\tilde{c}_n- \tilde{d}_n}) \sinh n\Gz \sinh \Gz \big)\cos n\Gt.
\end{align}
From the expressions of $c_n$ and $d_n$ given in Lemma \ref{lem:Phi_2}, one can see that,  for $n \geq 2$,
\begin{align*}
{\tilde{c}_n+\tilde{d}_n} &=\frac{4a}{s} \frac{ns e^{-ns}\cosh n s  - (ns)^2\eta_{2} + (ns)^3 \eta_{1}}{\sinh 2ns + 2ns\eta_{2} } - 4K_v  s \eta_{1} \frac{ n s }{\sinh 2ns + 2ns\eta_{2}},
\\
{\tilde{c}_n-\tilde{d}_n} &=\frac{4a}{s^2} \frac{(ns)^2 e^{-ns}\cosh n s   - (ns)^3\eta_{2}+  (ns)^2 s^2\eta_{1}}{\sinh 2ns +  2ns\eta_{2}} + 4K_v  \frac{ e^{-ns}\sinh ns +  n s\eta_{2}}{\sinh 2ns +  2ns\eta_{2}},
\end{align*}
where $\eta_1$ and $\eta_2$ are the quantities given in \eqref{eq:def_eta1_eta2}.
Define, for $0<x<\infty$,
\begin{align}
g_1(x) &:= \frac{x e^{-x}\cosh x  - x^2\eta_{2} + x^3 \eta_{1}}{\sinh 2x + 2x\eta_{2}}, \nonumber
\\
g_2(x) &:= \frac{x}{\sinh 2x + 2x\eta_{2}}, \nonumber
\\
g_3(x) &:=\frac{x^2 e^{-x}\cosh x   - x^3\eta_{2}+  x^2 s^2\eta_{1}}{\sinh 2x +  2x\eta_{2}}, \nonumber
\\
g_4(x) &:=\frac{e^{-x}\sinh x+ x \eta_2}{\sinh 2x + 2x\eta_{2}}. \label{gjdef}
\end{align}
Then, we have
\begin{align}
    {\tilde{c}_n+\tilde{d}_n} &=\frac{4a}{s} g_1(ns) - 4K_v  s \eta_{1} g_2(ns),\nonumber
\\
{\tilde{c}_n-\tilde{d}_n} &=\frac{4a}{s^2} g_3(ns) + 4K_v  g_4(ns).\label{an_bn_tilde_even_odd}
\end{align}

It follows from \eqref{eq:Ecal_GzGz_v2_symm} that
\begin{align}
     \Ecal_{\Gz\Gz}[\Bv_2-\BU] &=
    -K_v \frac{\sinh\Gz}{a}\sin \Gt +  h(\Gz,\Gt)2 c_1 \sinh 2\Gz \sin \Gt\nonumber
    \\
    &\quad + \frac{4a}{s}h(\Gz,\Gt)\sum_{n=2}^\infty \big(  g_1(ns) \sinh n \Gz \cosh \Gz \big)\sin n\Gt\nonumber
    \\&\quad
    - 4 K_v h(\Gz,\Gt)\sum_{n=2}^\infty \big(  g_2(ns) \sinh n \Gz \cosh \Gz \big)\sin n\Gt \nonumber
    \\
    &\quad + \frac{4a}{s^2}h(\Gz,\Gt)\sum_{n=2}^\infty \big(   g_3(ns) \cosh n\Gz \sinh \Gz \big)\sin n\Gt\nonumber
    \\
    &\quad
    + 4K_v h(\Gz,\Gt)\sum_{n=2}^\infty \big(   g_4(ns) \cosh n\Gz \sinh \Gz \big)\sin n\Gt.
\end{align}
If we define
\begin{align*}
&T_j^{++}(\Gz,\Gt) := a h\sum_{n=2}^\infty g_j(ns) \cosh n\Gz \cos n\Gt,
\\
&T_j^{+-}(\Gz,\Gt) :=  a h\sum_{n=2}^\infty g_j(ns) \cosh n\Gz \sin n\Gt,
\\
&T_j^{-+}(\Gz,\Gt) := a h \sum_{n=2}^\infty g_j(ns) \sinh n\Gz \cos n\Gt,
\\
&T_j^{--}(\Gz,\Gt) := a h \sum_{n=2}^\infty g_j(ns) \sinh n\Gz \sin n\Gt,
\end{align*}
then the component $\Ecal_{\Gz\Gz} $ can be rewritten as
\begin{align}
    \Ecal_{\Gz\Gz}[\Bv_2-\BU]
    &=T_0 + \frac{4}{s} \cosh \Gz T_1^{--} - \frac{4K_v s \eta_1}{a}\cosh \Gz T_2^{--} \nonumber
    \\&\quad + \frac{4}{s^2} \sinh \Gz T_3^{+-} +  \frac{4K_v}{a}  \sinh \Gz T_4^{+-},
    \label{eq:Ecal_GzGz_Tj_v2}
\end{align}
where
\begin{equation}\label{eq:T0_def}
T_0 =-K_v \frac{\sinh\Gz}{a}\sin\Gt +  h(\Gz,\Gt)2 c_1 \sinh 2\Gz \sin \Gt.
\end{equation}
Similarly, $\Ecal_{\Gz\Gt}$ can be written as
\begin{align}
    \Ecal_{\Gz\Gt}[\Bv_2-\BU]
    &=\widetilde{T}_0 + \frac{4}{s} \cosh \Gz T_1^{++} - \frac{4K_v s \eta_1}{a}\cosh \Gz T_2^{++} \nonumber
    \\&\quad + \frac{4}{s^2} \sinh \Gz T_3^{-+} +  \frac{4K_v}{a}  \sinh \Gz T_4^{-+},
    \label{eq:Ecal_GzGt_Tj_v2}
\end{align}
where
\begin{align}
\widetilde{T}_0 &= K_v \frac{\cosh2\Gz -2\cosh\Gz\cos\Gt+ \cos 2\Gt}{2a}  \nonumber
\\& \quad +  h(\Gz,\Gt) d_0\cosh\Gz + h(\Gz,\Gt) 2c_1\cosh 2\Gz \cos\Gt.
\label{eq:T0p_def}
\end{align}

The proof the following lemma is given in appendix \ref{sec:appendxE}.
\begin{lemma}\label{lem:gj_estim}  For $j=1,2,3,4$, we have
\beq\label{gjest}
|g_j(x)|+ |g_j'(x)| + |g_j''(x)|\lesssim (1+x^3)e^{-2x}, \quad 0<x<\infty.
\eeq
\end{lemma}

We omit the proof of the following lemma since it can be proved using Lemma \ref{lem:gj_estim} in the
same way as Lemma \ref{lem:pre_asymp_Sj}:

\begin{lemma}\label{lem:pre_asymp_Tj}
For $j=1,2,3,4$,
\begin{align}
T_{j}^{+ +} &=- \frac{ 1}{2} g_j(2s)\cosh 2\Gz \cos \Gt  + g_j(2s)\cosh 2\Gz \cos 2\Gt \nonumber \\
&\quad-\frac{1}{2} g_j(3s)\cosh 3\Gz \cos 2\Gt + O(s), \\
T_{j}^{+ -} &=- \frac{ 1}{2} g_j(2s)\cosh 2\Gz \sin \Gt  +g_j(2s)\cosh 2\Gz \sin 2\Gt \nonumber \\
&\quad  -\frac{1}{2} g_j(3s)\cosh 3\Gz \sin 2\Gt  + O(s), \\
T_{j}^{- \pm}&=O(s).	
\end{align}
\end{lemma}

We infer using the definition \eqnref{gjdef} of $g_j$ that for $k=2,3,$
\begin{align}
    &g_1(ks) = \frac{1}{4} + O(s),
\quad
    g_2(ks) = \frac{1}{4} + O(s), \nonumber
    \\
    &g_3(ks) = \frac{1}{4}s + O(s^2),
    \quad
    g_4(ks) = \frac{1}{2} + O(s). \label{eq:gj_zero_asymp}
\end{align}	
Thus, we have
$$
T_3^{+-} = O(s)
$$
and
\begin{align*}
T_4^{+-} &= - \frac{1}{4} \cosh 2\Gz \sin \Gt  + \frac{1}{2}\cosh 2\Gz \sin 2\Gt-\frac{1}{4} \cosh 3\Gz \sin 2\Gt  + O(s) \\
&= -  \frac{1}{4}(\sin \Gt - \sin 2\Gt)  + O(s),
\end{align*}
as $s \to 0$. It then follows from \eqref{eq:Ecal_GzGz_Tj_v2} that
\beq\label{Ezzest}
\Ecal_{\Gz\Gz}[\Bv_2-\BU] = T_0 + \frac{K_v \Gz}{a}(- \sin\Gt +  \sin 2\Gt ) +O(1).
\eeq
Since
\beq\label{c1est}
c_1 = a(-1+\coth 2s) + K_v \frac{1}{1+e^{2s}} =  \frac{1-s}{2}K_v + \frac{a}{2s} +O(s),
\eeq
it follows from \eqref{eq:T0_def} that
\begin{align}
T_0 &=-K_v \frac{\sinh\Gz}{a}\sin\Gt +  \frac{\cosh \Gz}{a} 2 c_1 \sinh 2\Gz \sin \Gt - \frac{\cos\Gt}{a}2 c_1 \sinh 2\Gz \sin \Gt\nonumber
\\
&= - K_v \frac{\Gz}{a} \sin\Gt + 2 K_v \frac{\Gz}{a}\sin\Gt - K_v  \frac{\Gz}{a} \sin 2\Gt + O(1)\nonumber
\\
&=   K_v \frac{\Gz}{a}\sin\Gt - K_v  \frac{\Gz}{a} \sin 2\Gt + O(1). \label{eq:T0_asymp}
\end{align}
This together with \eqnref{Ezzest} yields the desired estimate
\beq\label{Ezzest2}
\Ecal_{\Gz\Gz}[\Bv_2-\BU] = O(1).
\eeq

Likewise, we use \eqref{eq:Ecal_GzGt_Tj_v2}, Lemma \ref{lem:pre_asymp_Tj} and \eqnref{eq:gj_zero_asymp} to ensure
$$
\Ecal_{\Gz\Gt}[\Bv_2-\BU]
=\widetilde{T}_0 + (-\frac{1}{2s} + \frac{K_v s}{2a})\cos\Gt +  (+\frac{1}{2s} - \frac{K_v s}{2a})\cos 2\Gt + O(1).
$$
Then, using the estimate
\beq\label{d0est}
d_0 = \frac{a}{\sinh s \cosh s +s} -K_v \frac{\sinh^2 s}{s+\cosh s \sinh s} = \frac{a}{2s}- K_v \frac{s}{2} + O(s)
\eeq
in addition to \eqnref{c1est}, we obtain
\begin{align}
{T}_0' &= \frac{K_v}{2a} - \frac{K_v}{a}  \cos\Gt +\frac{1}{2}\frac{K_v}{a}\cos 2\Gt + \frac{1-\cos\Gt}{a}(\frac{a}{2s} + \frac{K_v s}{2})\nonumber
\\&\quad + \frac{1-\cos\Gt}{a} (K_v -s K_v + \frac{a}{s})\cos\Gt + O(1)
\nonumber
\\
&=\frac{K_v}{a} ( \frac{1}{2}) + \frac{1}{2s}-\frac{K_v s}{2a} - \frac{1}{2}(\frac{K_v}{a} - \frac{K_v s}{a} + \frac{1}{s})   + ( - \frac{1}{2s}+\frac{K_v s}{2a} - \frac{K_v s}{a} + \frac{1}{s}) \cos\Gt
\nonumber
\\
&\quad  + (\frac{1}{2}\frac{K_v}{a}- \frac{1}{2}(\frac{K_v}{a} - \frac{K_v s}{a} + \frac{1}{s}) ) \cos 2\Gt +O(1)
\nonumber
\\
&=   ( \frac{1}{2s} - \frac{K_v s}{2a})\cos\Gt + (\frac{K_v s }{2a} -\frac{1}{2s})\cos2\Gt +O(1),\label{eq:T0p_asymp}
\end{align}
where we have used the following identity for the second equality:
$$
(1-\cos\Gt)\cos\Gt = -\frac{1}{2} + \cos\Gt -\frac{1}{2}\cos2\Gt.
$$
Thus, we have
\begin{align*}
\Ecal_{\Gz\Gt}[\Bv_2-\BU] = O(1).
\end{align*}
So far we proved that
\begin{align*}
\| \Ecal[\Bv_2]  \|_\infty \lesssim 1.
\end{align*}

We now prove that the pressure $q_2$ is bounded independently of $\Gd$. It was shown in \eqnref{qtwo}
\begin{align}
    q_2 &=\frac{2\mu }{a} d_0 {w}_1 + \frac{2\mu}{a}\sum_{n=1}^\infty((n+1)c_n - (n-1)d_n ){w}_n \nonumber
    \\&\quad - \frac{2\mu}{a}\sum_{n=1}^\infty n c_n {w}_{n+1}
     + \frac{2\mu}{a}\sum_{n=1}^\infty n d_n{w}_{n-1} , \nonumber
\end{align}
where $w_n(\Gz,\Gt):=\sinh n\Gz \sin n \Gt$. It can be rewritten as
\begin{align}
        q_2 &=
        \frac{2\mu }{a} (-d_0 + 2(c_1 + d_2)) w_1 +\frac{2\mu}{a}\sum_{n=2}^\infty \frac{1}{n}(\tilde{c}_n - \tilde{d}_n - \tilde{c}_{n-1} + \tilde{d}_{n+1})w_n,
    \label{q2_series2}
\end{align}
where $\tilde{c}_n = n(n+1) c_n$ and $\tilde{d}_n = n(n-1) d_n$ for $n\geq 2$.

Note that
\begin{align*}
    \tilde{c}_n - \tilde{d}_n - \tilde{c}_{n-1} + \tilde{d}_{n+1} &=\frac{1}{2} (\tilde{c}_{n+1}+\tilde{d}_{n+1})  -\frac{1}{2}(\tilde{c}_{n-1}+\tilde{d}_{n-1})
    \\
    &\quad + \frac{1}{2}\Big(2(\tilde{c}_n - \tilde{d}_n)   - (\tilde{c}_{n-1}-\tilde{d}_{n-1}) -(\tilde{c}_{n+1}-\tilde{d}_{n+1})\Big).
\end{align*}
It then follows from \eqref{an_bn_tilde_even_odd} that
\begin{align}
    \tilde{c}_n - \tilde{d}_n - \tilde{c}_{n-1} + \tilde{d}_{n+1} &= \frac{2a}{s} ( g_1((n+1)s) -  g_1((n-1)s)) \nonumber
    \\
    &\quad -2K_v s^2 \eta_1 ( g_2((n+1)s)- g_2((n-1)s)) \nonumber
    \\
    &\quad   -2a (g_3((n+1)s)-2g_3(ns) + g_3((n-1)s)) \nonumber
    \\
    &\quad
    - 2K_v s^2 ( g_4((n+1)s) - 2g_4(ns) + g_4((n-1)s)). \label{+-+-cd}
\end{align}
By the mean value theorem, there are $x_{j,n} \in ((n-1)s, (n+1)s)$ such that
$$
|g_j((n+1)s) -  g_j((n-1)s)| \lesssim s |g_j'(x_{j,n})|, \quad j=1,2,
$$
and
$$
|g_j((n+1)s) - 2g_j(ns) + g_j((n-1)s))| \lesssim s^2 |g_j''(x_{j,n})|, \quad j=3,4.
$$
We then infer from \eqnref{gjest} that
$$
|g_j((n+1)s) -  g_j((n-1)s)| \lesssim s (1+(ns)^3)e^{-2ns}, \quad j=1,2,
$$
and
$$
|g_j((n+1)s) - 2g_j(ns) + g_j((n-1)s))| \lesssim s^2 (1+(ns)^3)e^{-2ns}, \quad j=3,4.
$$
Since $a\approx s$ and $K_v =O( s^{-1})$, it then follows from \eqnref{+-+-cd} that
$$
|\tilde{c}_n - \tilde{d}_n - \tilde{c}_{n-1} + \tilde{d}_{n+1}| \lesssim s (1+(ns)^3)e^{-2ns},
$$
and from \eqnref{q2_series2} that
\beq\label{2000}
|q_2| \lesssim (|d_0| + 2|c_1 + d_2|) + \sum_{n=2}^\infty \frac{1}{n}(1+(ns)^3)e^{-2ns} \sinh n \Gz.
\eeq

One can see from \eqnref{d0est} that
$$
d_0 = \frac{a}{2s}-\frac{K_v s}{2} +O(s) = O(1).
$$
One can also see from expressions of $c_1$ and $d_2$ in Lemma \ref{lem:Phi_2} that
$$
c_1 = \frac{1}{2}K_v +O(1), \quad d_2= - \frac{1}{2}K_v + O(1) ,
$$
and hence $c_1+d_2 = O(1)$. Thus we have from \eqnref{2000}
\begin{align*}
    |q_2| \lesssim 1 + \sum_{n=2}^\infty \frac{s\sinh n s}{n s}(1+(ns)^3)e^{-2ns} \lesssim 1+ \int_0^\infty \frac{\sinh t }{t}(1+t^3)e^{-2t} dt \lesssim 1.
\end{align*}
This completes the proof.
\qed

\section*{Concluding remarks}

In this paper, we have investigated the problem of quantifying the stress concentration in the narrow region between two rigid cylinders and derived precise estimates for the stress blow-up in the Stokes system when inclusions are circular cylinders of the same radii. We have shown that, even though the divergence of the velocity is confined to be zero, either the pressure component or the shear stress component of the stress tensor always blows up, and that the blow-up rate is $\Gd^{-1/2}$, where $\Gd$ is the distance between the cylinders. This blow-up rate coincides with the ones for elasto-statics and elasto-statics. In the course of deriving the results, it is proved that the blow-up of the stress tensor does not occur when the no-slip boundary is prescribed. We also derived an asymptotic decomposition formula which explicitly characterizes the singular behaviour of the solution. This formula may play an important role in computing the Stokes flow in presence of closely located rigid cylinders.

Since the method of bipolar coordinates is employed, extension of this paper's results to the case of circular cylinders with different radii is not a big issue. However, it is quite challenging to extend them to the more general case when the cross sections of the cylinders are strictly convex. In particular, proving no blow-up for the problem with the no-slip boundary condition on the convex boundaries seems already quite challenging.

\appendix

\section{Proof of Lemma \ref{lem:grad_eGz_eGt_estim}}\label{appendixC}

We have from \eqref{exey} that
$$
\Be_\Gz = \Ga \Be_x - \Gb  \Be_y, \quad \Be_\Gt = -\Gb \Be_x - \Ga  \Be_y.
$$
So we have
$$
|\nabla \Be_\Gz| + |\nabla \Be_\Gt| = 2( |\nabla \Ga| + |\nabla \Gb|) \lesssim |h \p_\Gz \Ga| + |h\p_\Gt\Ga| + |h\p_\Gz\Gb|	+|h\p_\Gt\Gb|	.
$$
Since
$$
h \p_\Gz \Ga = 	-\frac{\sinh\Gz \sin^2\Gt}{a(\cosh\Gz-\cos\Gt)},
$$
one can see that
$$
|h \p_\Gz \Ga| \lesssim |\Gt| \le |\Gz| + |\Gt|.
$$
Similarly one can show that
$$
|h\p_\Gt\Ga| , |h\p_\Gz\Gb|	, |h\p_\Gt\Gb| \lesssim |\Gz| + |\Gt|.
$$
This completes the proof. \qed

\section{Proof of Lemma \ref{lem:asymp_I1_J1}}\label{appendixA}

According to the transition relation \eqnref{Gsrel}, the stress tensor $\Gs[\Bh_1,p_1]$ is given by
$$
\Gs[\Bh_1,p_1] = \Xi \begin{bmatrix}
\Gs_{1,\Gz\Gz} &   \Gs_{1,\Gz\Gt}
\\
\Gs_{1,\Gz\Gt} & \Gs_{1,\Gt\Gt}
 \end{bmatrix} \Xi.
$$
In particular, we have
\beq
\Be_\Gz \cdot \Gs[\Bh_1,p_1] \Be_\Gz = \Gs_{1,\Gz\Gz}, \quad \Be_\Gz \cdot \Gs[\Bh_1,p_1] \Be_\Gt = \Gs_{1,\Gz\Gt}.
\eeq
On $\p D_2$ which is parametrized by $\{\Gz=s\}$, the outward unit normal $\nu$ is given by
$$
\nu |_{\p D_2} = - \Be_\Gz|_{\Gz=s},
$$
and, according to \eqnref{exey}, $\Be_x$ is expressed as $\Be_x=\Ga(\Gz,\Gt) \Be_{\Gz} -\Gb(\Gz,\Gt) \Be_\Gt$, where $\Ga$ and $\Gb$ are defined by \eqnref{pqdef}.
So, we have
\begin{align*}
\Be_x \cdot \Gs [\Bh_1,p_1]\nu = -\big( \Ga(s,\Gt) \Gs_{1,\Gz\Gz} -\Gb(s,\Gt) \Gs_{1,\Gz\Gt} \big).
\end{align*}
Due to \eqnref{lineele},
we have
\beq\label{A2}
\Ical_1 = - \int_{-\pi}^\pi \big( \Ga(s,\Gt) \Gs_{1,\Gz\Gz} -\Gb(s,\Gt) \Gs_{1,\Gz\Gt} \big) h(s,\Gt)^{-1} d\Gt.
\eeq

We now compute $( \Ga(s,\Gt) \Gs_{1,\Gz\Gz} -\Gb(s,\Gt) \Gs_{1,\Gz\Gt} ) h(s,\Gt)^{-1}$. It follows from the formulas \eqnref{eq:strain_bipolar1} and \eqnref{eq:strain_bipolar3} of the strain in bipolar coordinates, the strain-stress relation \eqref{eq:sigma_formula_bipolar}, and the formula \eqref{eq:stream_singular1} for the stream function that
\begin{align}
\Gs_{1,\Gz\Gz}|_{\Gz=s} &= \frac{A_2 \mu}{a} (2+\mbox{sech}\, 2s - 4 \cosh^3 s\, \mbox{sech}\, 2s \cos\Gt + \cos 2\Gt) , \label{4000}
\\
\Gs_{1,\Gz\Gt}|_{\Gz=s} &= -\frac{A_2\mu}{a}{(\cosh s-\cos\Gt)2\tanh 2s \sin\Gt}. \label{4010}
\end{align}
Using the definitions \eqnref{pqdef} of $\Ga$ and $\Gb$, \eqnref{4000} and \eqnref{4010}, we arrive at
\begin{align}
& \big( \Ga(s,\Gt) \Gs_{1,\Gz\Gz} -\Gb(s,\Gt) \Gs_{1,\Gz\Gt} \big) h(s,\Gt)^{-1}   = 2\mu A_1 (-1+\cosh s \,\mbox{sech}\,2s \cos\Gt). \label{Icalint}
\end{align}
Then by integrating both sides of \eqnref{Icalint} over $[-\pi,\pi]$, we obtain
\beq\label{Icalone2}
\Ical_1 = -4\pi \mu A_1 = -\frac{4\pi \mu}{2s-\tanh 2s} = -\frac{3\pi\mu}{2} \frac{1}{s^3} + O(s^{-1}).
\eeq
Thanks to the asymptotic formula \eqnref{sGd} of $s$ as $\Gd$ tends to $0$, we get the asymptotic formula \eqnref{Icalone} for $\Ical_1$.

Next we consider $\Jcal_1$. Similarly to the case of $\Ical_1$, we have
\begin{align*}
\int_{\p D_2} \BU \cdot \Gs[\Bh_1,p_1] \nu \, dl &
=\int_{\p D_2} \BU\cdot \Gs[\Bh_1,p_1]  (-\Be_\Gz) \, dl \\
&= - \int_{-\pi}^\pi (U_\Gz \Be_\Gz + U_\Gt \Be_\Gt)\cdot (\Gs_{1,\Gz\Gz} \Be_\Gz + \Gs_{1,\Gz\Gt}\Be_\Gt ) h(s,\Gt)^{-1} d\Gt \\
&=- \int_{-\pi}^\pi (U_\Gz \Gs_{1,\Gz\Gz} + U_\Gt \Gs_{1,\Gz\Gt})|_{\Gz=s} h(s, \Gt)^{-1} d\Gt .
\end{align*}
Since $U_\Gz=U \cdot \Be_\Gz$ and $U_\Gt=U \cdot \Be_\Gt$, it follows from \eqref{eq:bipolar_x_y} and \eqref{pqdef}  that
\begin{align}
U_\Gz|_{\Gz=s} &=\frac{a \sinh s \left(1-\cosh s \cos \Gt+\sin ^2\Gt\right)}{(\cosh s - \cos \Gt)^2},\nonumber \\
U_\Gt|_{\Gz=s} &=\frac{a \sin \Gt \left(1-\cosh s \cos \Gt-\sinh^2 s\right)}{(\cosh s- \cos \Gt)^2}.
\label{eq:UGzUGt_s}
\end{align}
It is convenient to use the following functions:
\beq\label{qnk}
q_n = q_n(s, \Gt) := \frac{\cos n\Gt}{\cosh s-\cos\Gt}, \quad n=0,1,2,\ldots.
\eeq
We obtain, by using \eqref{4000}, \eqref{4010} and \eqref{eq:UGzUGt_s}, that
\begin{align}
&(U_\Gz \Gs_{1,\Gz\Gz} + U_\Gt \Gs_{1,\Gz\Gt})h(s,\theta)^{-1} \nonumber \\
&=-a\mu A_1 \,\mbox{sech}\,2s \sinh s ((-1+2\cosh 2s)q_0 -2\cosh s q_1 + q_2) . \label{eq:USh}
\end{align}
As before, by integrating both sides of the equality in \eqref{eq:USh} over $[-\pi,\pi]$, we arrive at
\begin{align}
\mathcal{J}_1 &=
-a\mu A_1 \,\mbox{sech}\,2s \sinh s ((-1+2\cosh 2s)\Qcal_1^0 -2\cosh s \Qcal_1^1 + \Qcal_1^2), \label{J1_pre}
\end{align}
where
\beq\label{eq:def_Qnk}
\Qcal_{n}=\Qcal_n(s) := \int_{-\pi}^\pi q_n(s, \Gt) d\Gt.
\eeq

We also obtain the following asymptotic expansion of $\Qcal_n$, whose proof will be given after the current proof is completed.

\begin{lemma}\label{Qcalest}
It holds that
\beq\label{Qn}
\Qcal_n(s) =  \frac{2\pi}{s} - 2n\pi + (n^2-1/3)\pi s + O(s^2), \quad \mbox{as } s \to 0.
\eeq
\end{lemma}
Then, together with the asymptotic formula \eqnref{sGd} of $s$ as $\Gd$ tends to $0$, applying Lemma \ref{Qcalest} to \eqref{J1_pre} yields \eqnref{Jcalone}.
\qed

\medskip
\noindent{\sl Proof of Lemma \ref{Qcalest}}.
 One can easily see that $Q_n(s)$ is the real part of the following contour integral:
$$
-2 \int_C \frac{z^n}{z^2-2 z \cosh s +1} dz,
$$
where $C$ is the unit circle. Then the residue theorem yields
\beq\label{eq:Qcal_n}
\Qcal_n(s) = \frac{2\pi e^{-n s}}{\sinh s},
\eeq
from which \eqnref{Qn} follows.
\qed

\section{The asymptotics of the boundary integrals}\label{sec:appendixB}
Here we compute the asymptotics of the boundary integrals $\Ical_{22},\Ical_{23}, \Jcal_2, \Ical_{\mathrm{rot}}$, and $\Jcal_{\mathrm{rot}}$, and prove Lemma \ref{lem:asymp_I2_J2}.

\subsection{A lemma}

The following result will be used to prove Lemma \ref{lem:asymp_I2_J2}.

\begin{lemma} \label{lem:stream_integral}
Suppose that a solution $(\Bv,q)$ to the Stokes system on the exterior region $D^e$ satisfies $(\Bv,q)\in \Mcal$, and that its corresponding stream function $\Psi$ is given by
\begin{align}
(h\Psi)(\Gz,\Gt) &= K (\cosh\Gz - \cos \Gt)\ln (2\cosh\Gz-2\cos\Gt)
+a_0 \cosh \Gz + d_0 \Gz \sinh \Gz \nonumber
\\
&  \quad +\sum_{n=1}^\infty \big(a_n \cosh (n+1) \Gz + b_n \cosh (n-1)\Gz \big)\cos n\Gt.
\quad
\end{align}
Then we have the following formulas for the boundary integrals:
\begin{align}
&\int_{\p D_2} \Bpsi_2 \cdot \Gs[\Bv,q]\nu =  d_0 4\pi \mu,
\\
&\int_{\p D_2} \Bpsi_3 \cdot \Gs[\Bv,q]\nu = K 4\pi \mu a.
\end{align}
where $a$ is the number defined in \eqnref{adef}.
\end{lemma}

\proof
Let us write in terms of bipolar coordinates the stress tensor $\Gs[\Bv,q]$ and the strain tensor $\Ecal[\Bv,q]$ as
$$
\Gs[\Bv,q] = \Xi \begin{bmatrix}
\Gs_{\Gz\Gz} &   \Gs_{\Gz\Gt}
\\
\Gs_{\Gz\Gt} & \Gs_{\Gt\Gt}
 \end{bmatrix} \Xi,\quad \Ecal[\Bv,q] = \Xi \begin{bmatrix}
\Ecal_{\Gz\Gz} &   \Ecal_{\Gz\Gt}
\\
\Ecal_{\Gz\Gt} & \Ecal_{\Gt\Gt}
 \end{bmatrix} \Xi.
$$
One can show in the same way as of deriving \eqnref{A2} that
\begin{align}
\Kcal_2:=\int_{\p D_2} \Bpsi_2\cdot \Gs[\Bv,q]\nu = \int_{-\pi}^\pi \big ( \Gb(s,\Gt) \Gs_{\Gz\Gz} + \Ga(s,\Gt) \Gs_{\Gz\Gt} \big) h(s,\Gt)^{-1} d\Gt.
\end{align}
Using \eqref{eq:Be_Gt_formula}, one can also see that
\begin{align}
\Kcal_3:=\int_{\p D_2} \Bpsi_3\cdot \Gs[\Bv,q]\nu &= -\frac{a}{\tanh s}\int_{\p D_2} \Bpsi_2\cdot \Gs[\Bv,q] \nu -\frac{a}{\sinh s} \int_{\p D_2}\Be_\Gt \cdot \Gs[\Bv,q]\nu \nonumber\\
&=-\frac{a}{\tanh s}\Kcal_2+ \frac{a}{\sinh s} \int_{-\pi}^\pi \Gs_{\Gz\Gt}|_{\Gz=s} h(s,\Gt)^{-1}d\Gt. \label{eq:K3_identity}
\end{align}

We assume for a moment that the stream function $\Psi$ is given by
\beq\label{eq:Psi_temp1}
(h\Psi)(\Gz,\Gt) = K (\cosh\Gz - \cos \Gt)\ln (2\cosh\Gz-2\cos\Gt).
\eeq
Applying the formula \eqnref{eq:Laplacian_Psi} for the Laplacian in bipolar coordinates, we see that $\mu\GD\Psi=0$. Together with the relation \eqref{eq:pressure_bipolar} between the pressure and the stream function and the condition $q\rightarrow 0$ as $|\Bx|\rightarrow \infty$, this implies that the corresponding pressure $q=0$. Then, by \eqnref{eq:strain_bipolar1}-\eqnref{eq:strain_bipolar3} (the strain-stream function relation) and \eqref{eq:sigma_formula_bipolar} (the stress-strain relation),  we obtain
\begin{align}
\Gs_{\Gz\Gz}|_{\Gz=s} &= -K\frac{2\mu}{a}\sinh s\sin\Gt, \label{eq:Szz_temp1}
\\
\Gs_{\Gz\Gt}|_{\Gz=s} &=\frac{K\mu}{a}(\sinh^2 s -\sin^2\Gt + (\cosh\Gz-\cos\Gt)^2). \label{eq:Szt_temp1}
\end{align}
We also have
\begin{align*}
 \big ( \Gb(s,\Gt) \Gs_{\Gz\Gz} + \Ga(s,\Gt) \Gs_{\Gz\Gt} \big)|_{\Gz=s}h(s,\Gt)^{-1} = -2\mu \cosh s\cos\Gt.
\end{align*}
Then, by integrating over $[-\pi,\pi]$, we arrive at
$$
\Kcal_2 = 0.
$$

We now consider $\Kcal_3$. We see from \eqref{eq:K3_identity} and \eqref{eq:Szt_temp1}  that
\begin{align*}
    \Kcal_3&=0 +\frac{a}{\sinh s} {K\mu} \int_{-\pi}^\pi \Big(\frac{\sinh^2 s - \sin^2 \Gt}{\cosh s-\cos \Gt} +   \cosh s - \cos\Gt\Big)\,d\Gt
    \\
    &= \frac{a}{\sinh s} {K\mu}( 2\pi \sinh s -  2\pi e^{-s} +  2\pi \cosh s )
\\
&=K 4\pi\mu a,
\end{align*}
where we have used \eqref{eq:Qcal_n} for the second equality.
So far, we have computed $\Kcal_2$ and $\Kcal_3$ when the stream function $\Psi$ is given by \eqref{eq:Psi_temp1}.

Next we assume that $\Psi$ is given by
\beq\label{eq:Psi_temp2}
(h\Psi)(\Gz,\Gt)=d_0 \Gz \sinh \Gz+a_0 \cosh \Gz   +\sum_{n=1}^\infty \big(a_n \cosh (n+1) \Gz + b_n \cosh (n-1)\Gz \big)\cos n\Gt.
\eeq
By symmetry and from the fact that $\Bh_2|_{\p D_i} = (-1)^i\frac{1}{2}\Psi_2$, we have
\begin{align*}
\Kcal_2 &= \int_{\p D_1} \frac{-1}{2}\Psi_2\cdot \Gs[\Bv,q]\nu +  \int_{\p D_2} \frac{1}{2}\Psi_2\cdot \Gs[\Bv,q]\nu
\\
&=  \int_{\p D^e} \Bh_2\cdot \Gs[\Bv,q]\nu.
\end{align*}
Then, by \eqref{eq:int_parts_formula} (the divergence theorem), we have
$$
\Kcal_2 = -2\mu\int_{ D^e}  \Ecal[\Bh_2]:\Ecal[\Bv].
$$
Recall from \eqref{Ezz2}-\eqref{Ezt2} that
\begin{align}
\Ecal_{\Gz\Gz}[\Bh_2] = 0, \quad
\Ecal_{\Gt\Gt}[\Bh_2] = 0, \quad
\Ecal_{\Gz\Gt}[\Bh_2] = h(\Gz,\Gt) A_2\cosh \Gz.
\end{align}
By \eqref{eq:strain_bipolar3}, one can easily check that
\begin{align}
\Ecal_{\Gz\Gt}[\Bv] &=   h(\Gz,\Gt) d_0 \cosh \Gz \nonumber
\\
&\quad + h(\Gz,\Gt)\sum_{n=1}^\infty (n(n+1)a_n \cosh(n+1)\Gz + n(n-1)b_n \cosh(n-1)\Gz) \cos n\Gt. \label{eq:strain_temp2}
\end{align}
So we obtain
\begin{align*}
    \Kcal_2 &=-2\mu\int_{ D^e}  \Ecal_{\Gz\Gz}[\Bh_2]\Ecal_{\Gz\Gz}[\Bv] + 2\Ecal_{\Gz\Gt}[\Bh_2]\Ecal_{\Gz\Gt}[\Bv] + \Ecal_{\Gt\Gt}[\Bh_2]\Ecal_{\Gt\Gt}[\Bv]
    \\
    &= -2\mu \int_{-s}^s \int_{-\pi}^\pi  2\Ecal_{\Gz\Gt}[\Bh_2] \Ecal_{\Gz\Gt}[\Bv] \frac{1}{h(\Gz,\Gt)^2}  d\Gt d\Gz
    \\
    &=-2\mu (2\pi) \int_{-s}^s 2A_2 d_0 \cosh^2\Gz \,d\Gz  = -d_0 4\pi\mu  A_2  (2s+\sinh 2s) = d_0 4\pi \mu .
\end{align*}
We now compute $\Kcal_3$. We have from \eqref{eq:K3_identity} and \eqref{eq:strain_temp2} that
\begin{align}
\Kcal_3 &=-\frac{a}{\tanh s}\Kcal_2+ \frac{2\mu a}{\sinh s} \int_{-\pi}^\pi \Ecal_{\Gz\Gt}|_{\Gz=s} h(s,\Gt)^{-1}d\Gt
\\
&=- \frac{a}{\tanh s} (d_0 4\pi \mu) +  \frac{2\mu a}{\sinh s}  \int_{-\pi}^\pi d_0 \cosh s d\Gt\nonumber
\\
&= 0.
\end{align}
The proof is completed.
\qed

\subsection{Proof of Lemma \ref{lem:asymp_I2_J2}}

Now we are ready to compute the asymptotics of integrals  $\Ical_{22},\Ical_{23}, \Jcal_2, \Ical_{\mathrm{rot}}$, and $\Jcal_{\mathrm{rot}}$.

We first consider $\Ical_{22}$. We see from Theorem \ref{thm:boundedstress_rot} and Proposition \ref{lem:htwo} that
\beq\label{eq:estim_h2_h2t}
\|\Gs[\Bh_2 - \widetilde\Bh_2, p_2 - \widetilde{p}_2] \|_\infty \lesssim 1.
\eeq
So we get
$$
\Ical_{22} =  \int_{\p D_2} \Bpsi_2 \cdot \Gs[\widetilde{\Bh}_2,\widetilde{p}_2]\nu +O(1).
$$
Therefore, Lemma \ref{lem:stream_integral} with $\Psi=\widetilde{\Psi}_2$ yields
$$
\Ical_{22}  =A_2 4\pi \mu
+O(1).
$$
Hence, since $s \approx \sqrt{\Gd}$, the asymptotic formula \eqnref{Atwo} for $A_2$ yields \eqref{eq:I22}.

We now consider $\Ical_{23}$. Recall from \eqref{eq:I23_another} that
$$
\Ical_{23} = \int_{\p D_2} \Bpsi_2\cdot {\Gs[ \Bh_{\mathrm{rot}},p_{\mathrm{rot}}]}  \big|_+ \nu.
$$
Then, Lemma \ref{lem:stream_integral} with $\Psi=\Psi_{\mathrm{rot}}$ yields
$$
\Ical_{23} = -\frac{ K_{\mathrm{rot}} \sinh^2 s}{\sinh s \cosh s +s} 4\pi\mu.
$$
Similarly, we have
$$
\Ical_{\mathrm{rot}}  = K_{\mathrm{rot}} 4\pi\mu a.
$$
Since $s \approx \sqrt{\Gd}$, the asymptotic formula \eqref{eq:Krot_asymp} for $K_{\mathrm{rot}}$ yields \eqref{eq:I23} and \eqref{eq:Irot}.

Next we consider $\Jcal_2$ and $\Jcal_{\mathrm{rot}}$. Using the symmetry and the fact that $(\Bv_2-\BU)|_{\p D_e} = -\BU$, we have
\begin{align*}
    \Jcal_2 = \frac{1}{2}\int_{\p D^e}(-1) (\Bv_2-\BU) \cdot \Gs [\Bh_2,p_2]\nu .
\end{align*}
Here $\BU(x,y)=\BUsh=(y,x)^T$. Thanks to Green's formula \eqref{eq:int_parts_formula2}, the following holds:
\begin{align*}
\Jcal_2 &=
    -\frac{1}{2}\int_{\p D^e} \Bh_2 \cdot \Gs [\Bv_2-\BU,q_2]\nu  =  -\frac{1}{2}\int_{\p D_2} \Bpsi_2 \cdot \Gs [\Bv_2-\BU,q_2]\nu.
\end{align*}
It can be proved in the same way as the proof of Lemma \ref{lem:c21} that
$$
\int_{\p D^e} \Bh_2 \cdot \Gs [\Bv_2-\BU,q_2]\nu =0.
$$
Thus,
$$
\Jcal_2 = -\frac{1}{2}\int_{\p D_2} \Bpsi_2 \cdot \Gs [\Bv_2-\BU,q_2]\nu.
$$
Similarly, we have
$$
\Jcal_{\mathrm{rot}} = -\int_{\p D_2} \Bpsi_3 \cdot \Gs [\Bv_2-\BU,q_2] \nu .
$$

Since the stream function $\Psi_{v,2}$ associated with $(\Bv_2-\BU,q_2)$ is given in \eqnref{eq:W_biharmonic_solution}, we may apply  Lemma \ref{lem:stream_integral} to have
\begin{align*}
    \Jcal_2 &= -\frac{4\pi \mu}{2} \left(\frac{a}{\sinh s \cosh s +s} -K_v \frac{\sinh^2 s}{s+\cosh s \sinh s}\right),
\\
\Jcal_{\mathrm{rot}} &= - K_v 4\pi\mu a.
\end{align*}
Since $s \approx a \approx \sqrt{\Gd}$, the asymptotic formula \eqnref{eq:Kv_asymp} for $K_v$ yields \eqref{eq:J2} and \eqref{eq:Jrot}.
The proof is then completed.
\qed

\section{Proof of Lemma \ref{lem:fj_v1_estim}}\label{appendixD}

If $1<x<\infty$, then one can easily see that
$$
|f_j(x)|+ |f_j'(x)| + |f_j''(x)|\lesssim x^3e^{-2x}
$$
for $j=1,2$. So we consider the case when $2s \le x \le 1$, and prove
\begin{equation}\label{eq:fj_estim2}
|f_j(x)|+ |f_j'(x)| + |f_j''(x)|\lesssim 1, \quad j=1,2,
\end{equation}
and \eqnref{eq:Df1_DDf2_asymp}.

Let
\begin{align*}
\Ga_1 (x) &= x e^{-x}\sinh x  - x^2\eta_{2} + x^3 \eta_1  ,
\\	
\Ga_2 (x) &= x^2 e^{-x}\sinh x  -  x^3 \eta_2 + s^2 x^2  \eta_1 , \\
\Gb (x) &= \frac{1}{\sinh 2x -2x\eta_2},
\end{align*}
so that the following relations hold:
\begin{align*}
	f_1(x) = \Ga_1(x) \Gb(x) , \quad f_2(x) = \Ga_2(x) \Gb(x).
\end{align*}

One can see from the definition \eqnref{eq:def_eta1_eta2} of $\eta_j$ that
\beq\label{eta2}
\eta_1 = 1+ O(s^2), \quad \eta_2 = 1+\frac{2}{3}s^2 + R_1(s),
\eeq
where the remainder term $R_1(s)$ satisfies
\beq\label{Rone}
|R_1(s)| \le \frac{4}{15} s^4,
\eeq
provided that $s$ is sufficiently small.

Suppose $2s<x\leq 1$. Since
$$
\Ga_1(x) = (1-\eta_2) x^2 + (-1+\eta_1)x^3 + \frac{2}{3}x^4 + O(x^5),
$$
we have
\beq\label{Gaone}
\Ga_1 =	\frac{2}{3} a + O(x^5), \quad \Ga_1'(x) =	\frac{2}{3} a' + O(x^4), \quad \Ga_1''(x) = \frac{2}{3}	a'' + O(x^3),
\eeq
where
\beq\label{aone}
a(x):= x^4- s^2 x^2 .
\eeq
Likewise, since
$$
\Ga_2 (x) = (1- \eta_2) x^3 + \frac{2}{3} x^5 - (x^4-s^2x^2 \eta_1) + O(x^6) ,
$$
we have
\beq\label{Gatwo}
\Ga_2 =	\tilde{a} + O(x^6), \quad \Ga_2' =	\tilde{a}' + O(x^5), \quad \Ga_2'' =	\tilde{a}'' + O(x^4),
\eeq
where
\beq\label{tildea}
\tilde{a}(x):= \big( -1+ \frac{2}{3} x \big) a(x).
\eeq

Let
$$
w(x):= \sinh 2x -2x\eta_2.
$$
so that $\Gb(x)= w(x)^{-1}$. Note that
$$
\sinh 2x = 2x + \frac{4}{3} x^3 + R_2(x),
$$
where the remainder term $R_2$ satisfies $R_2(x) = O(x^5)$ and
\beq\label{Rtwo}
R_2 (x)\geq \frac{4}{15} x^5.
\eeq
Then
$$
w(x) = \frac{4}{3}x^3 - \frac{4}{3}s^2 x + R,
$$
where $R:= R_2(x)-R_1(s)x$. Since $x > 2s$, it follows from \eqnref{Rone} and \eqnref{Rtwo} that
$$
R \ge \frac{4}{15} x(x^4-s^s) \ge C x^5
$$
for some positive constant $C$. Therefore, we have
$$
w(x)= \frac{4}{3} b(x) (1+O(x^2)),
$$
where the remainder term $O(x^2)$ is larger than $C x^2$ for some positive constant $C$ and
\beq\label{b}
b(x)= x^3 - s^2 x = \frac{a(x)}{x}.
\eeq
Thus we have
\beq\label{Gbone}
\Gb(x)= w(x)^{-1} = \frac{3}{4} \frac{1}{b(x)} + O(x^{-1}).
\eeq
Since $\Gb'=-\Gb^2 w'$ and $\Gb''= 2\Gb^3 (w')^2 -\Gb^2 w''$, we have
\beq\label{Gbtwo}
\Gb'(x)= - \frac{3}{4} \frac{b'}{b^2} + O(x^{-2}), \quad \Gb''(x)= \frac{3}{4} \frac{2(b')^2-bb''}{b^3} + O(x^{-3}).
\eeq

Now it is easy to see that $f_1(x)= O(x)$ and $f_1'(x)=O(1)$. To prove the first part of \eqnref{eq:Df1_DDf2_asymp}, we invoke \eqnref{Gaone}, \eqnref{Gbone} and \eqnref{Gbtwo} to derive
$$
f_1' = \frac{1}{2} \frac{a'b-ab'}{b^2} + O(x).
$$
Since $a=xb$, we have
\beq\label{abab}
\frac{a'b-ab'}{b^2} =1,
\eeq
which yields the first part of \eqnref{eq:Df1_DDf2_asymp}. To prove that $f_1''$ is bounded, we again use \eqnref{Gaone}, \eqnref{Gbone} and \eqnref{Gbtwo} to derive
$$
f_1''= \frac{1}{2} \frac{a''b^2-2a'b'b + 2a(b')^2 - abb''}{b^3} + O(1).
$$
One can easily see that
$$
\frac{a''b^2-2a'b'b + 2a(b')^2 - abb''}{b^3} = \left( \frac{a'b-ab'}{b^2} \right)'.
$$
Thus, thanks to \eqnref{abab}, we infer
\beq\label{ababab}
\frac{a''b^2-2a'b'b + 2a(b')^2 - abb''}{b^3} =0.
\eeq
This proves \eqnref{eq:fj_estim2} for $j=1$.

It is easy to see that $f_2(x)= O(x)$ and $f_2'(x)=O(1)$. On the other hand, we have
$$
f_2'' = \frac{3}{4} \frac{\tilde{a}''b^2-2\tilde{a}'b'b + 2\tilde{a}(b')^2 - \tilde{a}bb''}{b^3} + O(x).
$$
Because of \eqnref{tildea}, \eqnref{abab} and \eqnref{ababab}, we have
$$
\frac{\tilde{a}''b^2-2\tilde{a}'b'b + 2\tilde{a}(b')^2 - \tilde{a}bb''}{b^3} =
\frac{\frac{4}{3} a'b^2 - \frac{4}{3} a b'b }{b^3} = \frac{4}{3}.
$$
Thus $f_2''=1+ O(x)$, which proves the second part of \eqnref{eq:Df1_DDf2_asymp} as well as \eqnref{eq:fj_estim2} for $j=2$. This completes the proof.
\qed

\section{Proof of Lemma \ref{lem:gj_estim}}\label{sec:appendxE}

The functions $g_j$ can be rewritten as
\begin{align*}
g_1(x) &=  \big(e^{-x} \cosh x - \eta_2 x  + \eta_1 x^2\big) v(x),
\\
g_2(x) &= v(x),
\\
g_3(x) &= \big(x e^{-x}\cosh x - x^2 \eta_2 + s^2 \eta_1 x\big)v(x),
\\
g_4(x) &=\Big(\frac{e^{-x}\sinh x}{x} + \eta_2\Big) v(x),
\end{align*}
where
\beq
v(x) = \frac{x}{\sinh 2x + 2x \eta_2}, \quad x>0.
\eeq

We estimate $v$ first. Since $\eta_2 = \sinh(2s)/(2s)\geq 1$, we have
$$
|v(x)| \leq  \frac{x}{\sinh 2x + 2x},
$$
and hence
\beq\label{vest}
|v(x)| \lesssim (1+x )e^{-2x}.
\eeq
By straight-forward computations, one can see that
$$
v'(x) = \Gg_1(x) (v(x))^2, \quad  v''(x)= \Gg_2(x) (v(x))^3,
$$
where
\begin{align*}
\Gg_1(x) &:= \frac{\sinh 2x - 2x \cosh 2x }{x^2},
\\
\Gg_2(x) &:= \frac{2(3x + x \cosh 4x - \sinh 4x) + 4\eta_2 ( 2x \cosh 2x  -(1+2x^2) \sinh 2x)}{x^3}.
\end{align*}

By Taylor expansions, it is easy to see that $\Gg_1$ and $\Gg_2$ are bounded if $0<x \le 1$. It is also easy to see that $|\Gg_1(x)| \lesssim x^{-1} e^{2x}$ and $|\Gg_2(x)| \lesssim x^{-2} e^{4x}$ if $1<x<\infty$. Putting these estimates together, we have
$$
|\Gg_1 (x)|\lesssim \frac{e^{2x}}{1+x}, \quad  |\Gg_2(x)| \lesssim  \frac{e^{4x}}{1+x^2}, \quad  0<x<\infty.
$$
Then, from \eqref{vest}, we obtain
\begin{align}
    |v'(x)| &\lesssim \frac{e^{2x}}{1+x} |v(x)|^2 \lesssim  (1+x)e^{-2x}, \label{v1est}
    \\
    |v''(x)| &\lesssim \frac{e^{4x}}{1+x^2} |v(x)|^3 \lesssim  (1+x) e^{-2x}. \label{v2est}
\end{align}

Since $g_2=v$, the estimate \eqnref{gjest} for $j=2$ is already proved. Let us prove it for $j=1$. We are ready to estimate $g_j$ and their derivatives. We consider $g_1$ only for simplicity. We write
$$
g_1(x) = \Gg(x) v (x), \quad \mbox{where}\quad  \Gg(x) = e^{-x}\cosh x - \eta_2 x + \eta_1 x^2.
$$
It is easy to show that
\beq\label{Ggest}
|\Gg(x)| \lesssim 1+x^2, \quad |\Gg'(x)|\lesssim 1+x, \quad |\Gg''(x)|\lesssim 1,
\eeq
and the estimate \eqnref{gjest} for $j=1$ is an easy consequence of \eqnref{vest}-\eqnref{Ggest}. \eqnref{gjest} for $j=3,4$ can be proved in the same way.
\qed


\begin{thebibliography}{99}


\bibitem{ACKLY-ARMA-13} H. Ammari, G. Ciraolo, H. Kang, H. Lee and K. Yun, Spectral analysis of the Neumann-Poincar\'e operator and characterization of the stress concentration in anti-plane elasticity, Arch. Rational Mech. Anal., 208 (2013), 275--304.

\bibitem{AKL} H. Ammari, H. Kang and M. Lim, Gradient estimates for solutions to the conductivity problem, Math. Ann. 332(2) (2005), 277--286.


\bibitem{BE-PF-1965} M. Bentwich and C. Elata, Eddy Formation in an Eccentric Annular Domain, Physics Fluids 8 (1965) 2204.

\bibitem{Batchelor-1967} G.K. Batchelor, {\sl An Introduction to Fluid Dynamics}, Cambridge University Press, 1967.


\bibitem{BLY-ARMA-09}  E.S. Bao, Y. Li and B. Yin,
Gradient estimates for the perfect conductivity problem, Arch. Rational Mech. Anal. 193 (2009), 195-226.

\bibitem{BLL-ARMA-15} J. Bao, H. Li and Y. Li, Gradient estimates for solutions of the Lam\'e
system with partially infinite coefficients, Arch. Rational Mech. Anal. 215 (2015), 307--351.

\bibitem{Berlyand-SIMA-06} L. Berlyand, L. Borcea, and A. Panchenko, Network approximation for effective viscosity of concentrated suspensions with complex geometry, SIAM J. Math. Anal. 36 (2006), 1580–-1628.


\bibitem{Crowdy-IJNLM-2011}  D. G. Crowdy, Treadmilling swimmers near a no-slip wall at low Reynolds number, Intl J.
Non-Linear Mech. 46 (2011), 577–-585.


\bibitem{DRM-CJCE-1967}
C. L. Darabaner,  J. K. Raasch and  S. G. Mason, Particle motions in sheared suspensions XX: Circular cylinders, The Canadian Journal of Chemical Engineering, 45 (1967), 3--12.

\bibitem{FA-CES-1967}
N. A. Frankel and A. Akrivos, On the viscosity of a concentrated suspension of solid spheres. Chemical Engineering
Science, 22 (1967), 847--853.

\bibitem{Graham-ASR-1981} A.L. Graham, On the viscosity of suspension of solid spheres, Applied Scientific Research, 37 (1981), 275--286.


\bibitem{LLBY-QAM-14} H. Li, Y. Li, E.S. Bao and B. Yin,
Derivative estimates of solutions of elliptic systems in narrow domains, Q. Appl. Math. 72 (2014), 589--596.

\bibitem{IC-JFM-2017} K. Ishimoto and D. G. Crowdy,
Dynamics of a treadmilling microswimmer near
a no-slip wall in simple shear, J. Fluid Mech. 821 (2017), 647–-667.

\bibitem{Milton-book-2001} G.W. Milton, The Theory of Composites, Cambridge Monographs on Applied and Computational Mathematics, Cambridge University Press, 2001.



\bibitem{Jeffrey-PTRS-1921} G. B. Jeffery, Plane stress and plane strain in bipolar coordinates, Phil. Trans. Roy. Soc. London A 221 (1921), 265-293.

\bibitem{Jeffrey-PRSA-1922} G. B. Jeffrey, The rotation of two circular cylinders in a viscous fluid, Proc. Roy. Soc. A 101 (1922), 169--174.


\bibitem{JO-QJMAM-1981} D. J. Jeffrey and Y. Onishi, The slow motion of a cylinder next to a plane wall. The Quarterly Journal of Mechanics and Applied Mathematics, 34 (1981) 129–137.

\bibitem{Ju-QAM-19}
H. Ju, H. Li and L. Xu, Estimates for elliptic systems in a narrow region arising from composite materials, Q. Appl. Math. 77 (2019), 177-199.

\bibitem{NK-JFM-1984} K. C. Nunan and J. B. Keller, Effective viscosity of a periodic suspension, Journ. of Fluid Mech. 142 (1984) 269--287.

\bibitem{KY19} H. Kang and S. Yu, Quantitative characterization of stress concentration in the presence of closely spaced hard inclusions in two-dimensional linear elasticity, Arch. Rational Mech. Anal. 232 (2019), 121--196.

\bibitem{Mitrea-book-2012}
M. Mitrea and M. Wright, Boundary value problems for the Stokes system in arbitrary Lipschitz domains, Soci\'et\'e math\'ematique de France, 2012.


\bibitem{Raasch-PhD-1961} J. K. Raasch,  Beanspruchung und Verhalten suspendierter Feststoffteilchen in Scherstr\"omungen hoher Za\"ahigkeit, Ph.D. Thesis, Faculty of Mechanical Engineering, Karlsruhe Technical University, Karlsruhe, Germany 1961.

 \bibitem{Raasch-ZAMM-1961}  J. K. Raasch, Das Verhalten suspendierter Feststoffteilchen in Scherstr\"omungen hoher Z\"ahigkeit, Z. Angew. Math. Mech. 41 (1961) 147–-151.


\bibitem{Schubert-JFM-1967} G. Schubert, Viscous flow near a cusped corner, J. Fluid Mech. 27 (1967) 647-656.


\bibitem{Smythe-1968} W. R. Smythe, {\sl Static and dynamic electricity}, McGraw-Hill, New York, 1968.


\bibitem{Smith-M-1991} S. H. Smith, The rotation of two circular cylinders in a viscous fluid, Mathematika, 38 (1991) 63--66.


\bibitem{Wannier-QAM-1950} G. H. Wannier,  Hydrodynamics of lubrication, Q. Appl. Math. 8 (1950), 7-32.


\bibitem{Wakiya-JPSJ-1975} S. Wakiya, Application of bipolar coordinates to the two-dimensional creeping motion of a liquid. I. Flow over a projection or a depression on a wall. Journal of the Physical Society of Japan, 39 (1975), 1113-1120.

\bibitem{Wakiya-JPSJ-1975-II} S. Wakiya, Application of bipolar coordinates to the two-dimensional creeping motion of a liquid. II. Some problems for two circular cylinders in viscous fluid. Journal of the Physical Society of Japan, 39 (1975), 1603-1607.


\bibitem{Watson-M-1995}E. J. Watson, The rotation of two circular cylinders in a viscous fluid, Mathematika 42 (1995) 105--126.


\bibitem{Yun-SIAP-07} K. Yun, Estimates for electric fields blown up between closely adjacent conductors with arbitrary shape, SIAM J. Appl. Math. 67 (2007), 714--730.



\end{thebibliography}
\end{document}